\documentclass[a4paper]{amsart}
\usepackage[T1,T2A]{fontenc}
\usepackage[utf8]{inputenc}

\usepackage[sorting=nyt,giveninits=true]{biblatex} 
    \AtBeginBibliography{\small}
\usepackage{algorithm}
\usepackage{algorithmic}

\usepackage[colorlinks=True, citecolor=black, linkcolor=black, urlcolor=black]{hyperref}
\usepackage{caption}
\usepackage[nameinlink]{cleveref}
\usepackage{amsmath}
\usepackage{amssymb}
\usepackage{mathtools}
\usepackage{tikz}
\usepackage{tikz-cd}
\usepackage{booktabs}
\usepackage{comment}
\usepackage[shortlabels]{enumitem}

\DeclareCaptionType{algType}[Algorithm][List of mytype]
\crefname{algType}{algorithm}{algorithms}
\newenvironment{myalgorithm}{\captionsetup{type=algType}}{}

\theoremstyle{definition}
\newtheorem{definition}{Definition}[section]
\newtheorem{example}[definition]{Example}
\newtheorem{remark}[definition]{Remark}

\theoremstyle{plain}
\newtheorem{theorem}[definition]{Theorem}
\newtheorem{lemma}[definition]{Lemma}
\newtheorem{proposition}[definition]{Proposition}
\newtheorem{corollary}[definition]{Corollary}
\newtheorem{addendum}[definition]{Addendum}

\DeclareMathOperator{\Hom}{Hom}

\DeclareMathOperator{\Aut}{Aut}

\DeclareMathOperator{\sign}{sign}
\DeclareMathOperator{\discr}{discr}

\DeclareMathOperator{\tf}{tf}

\DeclareMathOperator{\rank}{rank}
\DeclareMathOperator{\std}{std}

\newcommand{\QQ}{\mathbb{Q}}

\newcommand{\ZZ}{\mathbb{Z}}

\newcommand{\IQ}{\mathbb{Q}}
\newcommand{\IR}{\mathbb{R}}
\newcommand{\IC}{\mathbb{C}}
\newcommand{\IZ}{\mathbb{Z}}
\newcommand{\IN}{\mathbb{N}}

\newcommand{\bA}{\mathbf{A}}
\newcommand{\bD}{\mathbf{D}}
\newcommand{\bE}{\mathbf{E}}
\newcommand{\bU}{\mathbf{U}}

\newcommand{\bu}{\mathbf{u}}
\newcommand{\bv}{\mathbf{v}}
\newcommand{\bw}{\mathbf{w}}
\newcommand{\bLambda}{\mathbf{\Lambda}}

\newcommand{\cK}{\mathcal{K}}
\newcommand{\cO}{\mathcal{O}}

\newcommand{\I}{\mathrm{I}}
\newcommand{\II}{\mathrm{II}}

\makeatletter
\newcommand*{\defeq}{\mathrel{\rlap{%
                     \raisebox{0.3ex}{\(\m@th\cdot\)}}%
                     \raisebox{-0.3ex}{\(\m@th\cdot\)}}%
                     =}
\makeatother

\title[Idoneal genera and K3 surfaces covering an Enriques surface]{Idoneal genera and K3~surfaces \\ covering an Enriques surface}
\thanks{The authors acknowledge the financial support of the following research projects: DFG SFB-TRR~45, DFG SFB-TRR~195, TÜBİTAK 118F413.}

\newcommand{\Lminus}{\hyperlink{Lminus}{\bLambda^-}}
\newcommand{\Etilde}{\hyperlink{tildeE}{\mathbf{\tilde{E}}}}

\author{Simon Brandhorst}
\address[Simon Brandhorst]{Fachrichtung Mathematik,
Universität des Saarlandes,
Campus E24, 66123 Saarbrücken, Germany}
\email{brandhorst@math.uni-sb.de}

\author{Serkan Sonel}
\address[Serkan Sonel]{Matematik Bolümü, Bilkent Üniversitesi, 06800 Ankara, Turkey}
\email{serkan.sonel@bilkent.edu.tr}

\author{Davide Cesare Veniani}
\address[Davide Cesare Veniani]{Institut für Diskrete Strukturen und Symbolisches Rechnen, Universität Stuttgart, Pfaffenwaldring 57, 70569 Stuttgart, Germany}
\email{davide.veniani@mathematik.uni-stuttgart.de}


\date{\today}

\makeatletter
\@namedef{subjclassname@2020}{
  \textup{2020} Mathematics Subject Classification}
\makeatother


\begin{document}

\begin{abstract}
    We introduce the notion of idoneal genera, which are a generalization of Euler's idoneal numbers. We prove that there exist only a finite number of idoneal genera, and we provide an algorithm to enumerate all idoneal genera of rank at least 3. As an application, we classify transcendental lattices of K3 surfaces covering an Enriques surface.
\end{abstract}

\maketitle

\section{Introduction}
In this paper, an (integral) \emph{lattice} of \emph{rank} (or \emph{dimension}) \(r\) is defined as a finitely generated free \(\IZ\)-module \(L \cong \IZ^r\) equipped with a nondegenerate symmetric bilinear pairing \(b \colon L \times L \to \IZ\).  
A \emph{genus} is a complete set of isomorphism classes of integral lattices that are equivalent over \(\IR\) and over \(\IZ_p\) for each prime \(p\) to a given lattice.  
The \emph{parity}, \emph{rank}, \emph{signature}, and \emph{determinant} of a genus \(g\) refer to the respective properties of any lattice \(L\) within \(g\).  
We denote by \([n]\) the lattice of rank \(1\) with a generator \(e\) such that \(b(e,e) = n\).

An \emph{idoneal} (or \emph{suitable} or \emph{convenient}) number is a positive integer \(n\) satisfying the following property: every odd number \(m\) coprime to \(n\) is a prime whenever \(m\) can be expressed as \(x^2 + n y^2\), where \(x\) and \(y\) are relatively prime, and the equation \(m = x^2 + ny^2\) has exactly one solution with \(x, y \geq 0\). This terminology originates from Euler \cite{euler:numeri.idonei}.  

We extend this concept to define idoneal genera as follows:

\begin{definition} \label{def:idoneal}
    A positive definite genus \(g\) is called \emph{idoneal} if every lattice \(L\) in \(g\) represents~\(1\), that is, contains an element of square \(1\) or, equivalently, is of the form \(L \cong L' \oplus [1]\) for some lattice \(L'\).
\end{definition}

\begin{remark} \label{rmk:idoneal}
It is known (see, e.g., \cite[Theorem 6]{Kani:idoneal.numbers}) that a positive integer \(n\) is an idoneal number if and only if the genus of \([n] \oplus [1]\) consists of a single class.
If a lattice \(L\) belongs to an idoneal genus \(g\) of rank \(2\) and determinant \(n\), then \(L \cong [n] \oplus [1]\). In particular, \(L\) is unique in \(g\).  
Therefore, a genus \(g\) of rank \(2\) is idoneal if and only if \(g = \{ [n] \oplus [1] \}\) for some idoneal number \(n\).  

This observation establishes a bijection between idoneal numbers and idoneal genera of rank \(2\). In this sense, the concept of an idoneal genus generalizes the notion of an idoneal number.
\end{remark}

\begin{example}
    The genus of the lattice \([11] \oplus [1]\) is not idoneal because it also contains the lattice with Gram matrix
    \[
        \begin{pmatrix} 4 & 1 \\ 1 & 3 \end{pmatrix},
    \]
    which does not represent \(1\). Indeed, \(11\) is the smallest non-idoneal number.
\end{example}

The problem of enumerating all idoneal genera was posed by Kani~\cite[Problem~45]{Kani:idoneal.numbers}.
The main theorem of the first part of the present paper is the following. 

\begin{theorem}[see §\ref{subsec:finiteness}] \label{thm:idoneal}
    There exist only finitely many idoneal genera.
\end{theorem}

\begin{remark} \label{rmk:GRH}
    By \Cref{rmk:idoneal}, the problem of enumerating idoneal genera of rank~\(2\) is equivalent to the problem of enumerating idoneal numbers. 
    There are exactly \(65\) idoneal numbers known, the largest one being \(1848\) (sequence A000926 in the OEIS \cite{OEIS:A000926}).
    Weinberger \cite{weinberger} proved that the list is complete if the generalized Riemann hypothesis holds.
    If it does not hold, then there could exist one or two more idoneal numbers (see \cite[Corollary~23]{Kani:idoneal.numbers}), hence one or two more idoneal genera of rank~\(2\).
\end{remark}

In this paper, we enumerate all idoneal genera of rank at least~\(3\) without the assumption of the generalized Riemann hypothesis. The enumeration is carried out using the computer algebra systems \texttt{Magma}~\cite{magma}, \texttt{PARI}~\cite{pari} and \texttt{sageMath}~\cite{sage}.
The list of all known \(577\) idoneal genera is contained in the ancillary file \verb+idoneal.genera.txt+ published on Zenodo~\cite{brandhorst_2024_10617125}.

\begin{addendum}[see §\ref{subsec:enumeration}] \label{add:idoneal}
    There are no idoneal genera of rank \(r > 13\). For each \(3 \leq r \leq 13\), there exist exactly \(|I_r|\) idoneal genera of rank~\(r\), with \(|I_r|\) given in \Cref{tab:summary}.
\end{addendum}

We also present an application of the notion of idoneal genus to the problem of characterizing complex K3 surfaces covering an Enriques surface. 

Recall that a smooth proper algebraic surface \(X\) defined over \(\IC\) such that \(H^1(X,\cO) = 0\) is called a \emph{K3 surface} if its canonical bundle \(\cK\) is trivial and it is called an \emph{Enriques surface} if \(\cK\) is not trivial, but \(\cK^{\otimes 2}\) is. We say that a K3 surface \(X\) \emph{covers} an Enriques surface \(Y\) if there exists a finite étale morphism \(X \to Y\) of degree~\(2\). 
The cohomology group \(H^2(X,\IZ)\) of a K3 surface~\(X\), together with the Poincaré pairing, is a unimodular lattice of rank~\(22\). 
It contains the \emph{Néron--Severi lattice}, defined as the image~\(S\) of the map \(H^1(X, \cO^*) \rightarrow H^2(X,\IZ)\) coming from the exponential sheaf sequence. The \emph{transcendental lattice} of \(X\) is the orthogonal complement of \(S\) in \(H^2(X,\IZ)\), and it has signature \((2, \lambda-2)\).

If \(e_1,\ldots,e_r \in L\) is a system of generators of a lattice $L$, the associated \emph{Gram matrix} is the symmetric matrix with entries \(b_{ij} = b(e_i,e_j)\).
The \emph{determinant} \(\det(L)\) is the determinant of any such matrix. A lattice~\(L\) is \emph{unimodular} if \(|{\det(L)}| = 1\).
A lattice~\(L\) is \emph{even} if \(e^2 = b(e,e) \in 2\IZ\) for each \(e \in L\), otherwise it is \emph{odd}.
An embedding of lattices \(L \hookrightarrow M\) is called \emph{primitive} if \(M/L\) is free; moreover, \(M\) is an \emph{overlattice} of~\(L\) of index~\(m\) if \(\rank M = \rank L\) and \(m = |M/L|\).
We write \(L(n)\) for the lattice with the pairing defined by the composition \(L \times L \rightarrow \IZ \xrightarrow{n} \IZ\) and we put
\[
n L \coloneqq L \oplus \ldots \oplus L \quad \text{(\(n\) times)}.
\]
We denote by \(\bU\) the unique unimodular, even, indefinite lattice of rank~\(2\), and by \(\bA_n,\bD_n,\bE_n\) the positive definite ADE lattices.
We set 
\[
    \hypertarget{Lminus}{\bLambda^-} \coloneqq \bU \oplus \bU(2) \oplus \bE_8(-2).
\]
(Note that \(\bE_8\) denotes a negative definite lattice in \cite{Keum90}.)

Our starting point is the following criterion proved by Keum (under an additional assumption which is actually superfluous, see~\cite{Ohashi}). 

\begin{theorem}[Keum's criterion \cite{Keum90}] \label{thm:keum}
A K3 surface \(X\) with transcendental lattice~\(T\) covers an Enriques surface if and only if 
there exists a primitive embedding \(T \hookrightarrow \Lminus\) 
such that there exists no vector \(v \in T^\perp\) with \(v^2 = -2\). \qed
\end{theorem}

Our aim is to restate Keum's criterion in a way that makes the condition on the transcendental lattice~\(T\) more explicit and easy to be checked, completing the work started by Sertöz~\cite{Sertoez05} and Ohashi~\cite{Ohashi}, and partially continued by Lee~\cite{Lee12} and Yörük~\cite{Yoruk}.

\begin{definition} \label{def:co-idoneal}
    A lattice \(T\) of signature \((2,\lambda-2)\) is called \emph{co-idoneal} if \(T\) embeds primitively into \(\Lminus\), and for each primitive embedding \(T \hookrightarrow \Lminus\) there exists \(v \in T^\perp\) with \(v^2 = -2\). 
\end{definition}

The name `co-idoneal' is justified by \Cref{prop:co-idoneal.lattices/idoneal.genera}.
According to Keum's criterion, there are two possible reasons why a K3 surface~\(X\) may fail to cover any Enriques surface: either \(T\) does not embed primitively into~\(\Lminus\), or \(T\) is a co-idoneal lattice.  
Using Nikulin's theory of discriminant forms, we establish the following theorem.

\begin{theorem}[see §\ref{subsec:proof-transc-thm}] \label{thm:transcendental}
Let \(\Etilde\) be a positive definite lattice of rank \(8\) and discriminant form~\(3\bu_1\), where \(\bu_1\) is the discriminant form of \(\bU(2)\).
If \(X\) is a K3 surface with transcendental lattice~\(T\) of rank~\(\lambda\), then \(X\) covers an Enriques surface if and only if \(T\) is not a co-idoneal lattice and
one of the following conditions holds: 
\begin{enumerate}[(i)]
    \item \label{cond:2<=lambda<=6} \(2 \leq \lambda \leq 6\) and \(T\) admits a Gram matrix of the form
    \[
    \begin{pmatrix} 
        2a_{11} & a_{12} & \ldots & a_{1\lambda} \\
        a_{12} & 2a_{22} & & \vdots \\
        \vdots &  & \ddots & \vdots \\
        a_{1\lambda} & \ldots & \ldots & 2a_{\lambda\lambda}
    \end{pmatrix}
    \]
    such that \(a_{ij}\) is even for each \(2 \leq i, j \leq \lambda\),
    \item \(\lambda = 7\) and there exists an even lattice \(T'\) with \(T \cong \bU \oplus T'(2)\),
    \item \(\lambda = 7\) and there exists a lattice \(T'\) with \(T \cong \bU(2) \oplus T'(2)\),
    \item \(\lambda = 8\) and there exists an even lattice \(T'\) with \(T \cong \bU \oplus \bU(2) \oplus T'(2)\),
    \item \(\lambda = 8\) and there exists a lattice \(T'\) with \(T \cong 2\bU(2) \oplus T'(2)\),
    \item \label{cond:lambda=9a} \(\lambda = 9\) and there exists an even lattice \(T'\) with \(\bU(2) \oplus T \cong \Etilde(-1) \oplus T'(2)\),
    \item \(\lambda = 9\) and there exists a lattice \(T'\) with \(\bU \oplus T \cong \Etilde(-1) \oplus T'(2)\).
    \item \label{cond:lambda=10.1} \(\lambda = 10\) and there exists an even lattice \(T'\) with \(T \cong \Etilde(-1) \oplus T'(2)\),
    \item \label{cond:lambda=10.2} \(\lambda = 10\) and there exists a lattice \(T'\) with \(T \cong \bE_8(-2) \oplus T'(2)\),   
    \item \(\lambda = 11\) and there exists \(n > 0\) with \(T \cong \bU \oplus \bE_8(-2) \oplus [4n]\),
    \item \(\lambda = 11\) and there exists \(n > 0\) with \(T \cong \bU(2) \oplus \bE_8(-2) \oplus [2n]\),
    \item \label{cond:lambda=12} \(\lambda = 12\) and \(T \cong \Lminus\).
\end{enumerate}
\end{theorem}
By \cite[Theorem~1.14.2]{Nikulin:int.sym.bilinear.forms}, the given conditions do not depend on the choice of the lattice \(\Etilde\) with properties as above.
One example of a lattice \(\hypertarget{tildeE}{\mathbf{\tilde E}}\) with the required properties is the lattice with the following Gram matrix:
{\small \[
\begin{pmatrix}
    2 & 0 & 0 & 0 & 0 & 1 &  1 &  1 \\
    0 & 2 & 0 & 0 & 0 & 1 &  1 &  1 \\
    0 & 0 & 2 & 0 & 0 & 1 &  1 &  1 \\
    0 & 0 & 0 & 2 & 0 & 1 &  1 &  1 \\
    0 & 0 & 0 & 0 & 2 & 1 & -1 & -1 \\
    1 & 1 & 1 & 1 & 1 & 4 &  1 &  1 \\
    1 & 1 & 1 & 1 & -1 & 1 & 4 &  2 \\
    1 & 1 & 1 & 1 & -1 & 1 & 2 & 4
\end{pmatrix}.
\]}

A characterization analogous to \Cref{thm:transcendental} was obtained by Morrison~\cite{morrison:K3.large.picard} for Kummer surfaces. 
Morrison restated Nikulin's criterion~\cite{nikulin:kummer} that a K3 surface~\(X\) with transcendental lattice~\(T\) is a Kummer surface if and only if there exists a primitive embedding \(T \hookrightarrow 3 \bU(2)\) (see \cite[Corollary 4.4]{morrison:K3.large.picard}).

According to Ohashi~\cite{Ohashi07}, the number of isomorphism classes of Enriques surfaces covered by a fixed K3 surface is finite. \Cref{thm:transcendental} classifies K3 surfaces for which this number is nonzero. For the related problem of computing this number explicitly, see \cite{Ohashi, Ohashi07, ShimadaVeniani}.
For an application of our results to the case of one-dimensional families of K3 surfaces, see \cite{Festi.Veniani}.

The last core result of this paper is the enumeration of all co-idoneal lattices. 
The key observation is the connection between idoneal genera and co-idoneal lattices contained in \Cref{prop:co-idoneal.lattices/idoneal.genera}:
the orthogonal complement in \(\Lminus\) of a co-idoneal lattice is always of the form \(L(-2)\), 
where \(L\) is a lattice belonging to a uniquely determined idoneal genus. 

\begin{theorem} \label{thm:co-idoneal}
    There exist only finitely many co-idoneal lattices.
\end{theorem}

\begin{addendum}[see §\ref{subsec:co-idoneal}] \label{add:co-idoneal}
For each \(\lambda \in \IZ_{\geq 2}\), there exist \(|E_\lambda|\) co-idoneal lattices of signature \((2, \lambda - 2)\), with \(|E_\lambda|\) given in \Cref{tab:summary} if \(\lambda \leq 11\) and \(|E_\lambda| = 0\) otherwise.
There exist no other co-idoneal lattices of rank \(\lambda \neq 10\), and there exist at most two more of rank~\(10\) if the generalized Riemann hypothesis does not hold.
\end{addendum}

Note that a co-idoneal lattice is always of the form \(T = T'(2)\) for some odd lattice~\(T'\), called the \emph{half} of \(T\) (\Cref{cor:half.exc}).
The list of all halves of the known \(550\) co-idoneal lattices is contained in the ancillary file \verb+half.co-idoneal.lat.txt+ published on Zenodo~\cite{brandhorst_2024_10617125}.

\begin{table}
\caption{The number \(|I_r|\) of idoneal genera of rank \(r\) and the corresponding number \(|E_\lambda|\) of co-idoneal lattices of rank~\(\lambda\) (see \Cref{add:idoneal}, \Cref{add:co-idoneal}, and \Cref{rmk:GRH} for *).}
\label{tab:summary}
\begin{tabular}{llllllllllllll}
\toprule
 \(r\) & 1 & 2 & 3 & 4 & 5 & 6 & 7 & 8 & 9 & 10 & 11 & 12 & 13 \\
 \(|I_r|\) & 1 & 65* & 110 & 122 & 107 & 76 & 47 & 24 & 13 & 6 & 4 & 1 & 1 \\
\midrule
 \(\lambda\) & 11 & 10 & 9 & 8 & 7 & 6 & 5 & 4 & 3 & 2 & -- & -- & -- \\
 \(|E_\lambda|\) & 1 & 65* & 110 & 122 & 107 & 76 & 41 & 17 & 8 & 3 & -- & -- & -- \\
\bottomrule
\end{tabular}
\end{table}

\subsection{Contents of the paper}
The paper is divided into two sections.
In \Cref{sec:idoneal.genera} the relevant facts about the Smith--Minkowski--Siegel mass formula are summarized. The section is then devoted to the classification of idoneal genera and contains the proof of \Cref{thm:idoneal}. 
In \Cref{sec:K3-covering-Enriques}, after recalling some results on finite discriminant forms, we first determine which transcendental lattices embed into~\(\Lminus\), proving \Cref{thm:transcendental}.
We conclude the paper with the enumeration of all co-idoneal lattices.

\subsection*{Acknowledgments}
The authors would like to warmly thank Alex Degtyarev, Markus Kirschmer, Stéphane Louboutin, Rainer Schulze-Pillot, Ali Sinan Sertöz and John Voight for sharing their insights.
The authors are also very grateful to the anonymous referees for their careful proofreading and useful suggestions.

\section{Idoneal genera}\label{sec:idoneal.genera}

This section aims to prove \Cref{thm:idoneal} and present an algorithm for enumerating all idoneal genera (\Cref{def:idoneal}) of rank \(n \geq 3\).

In §\ref{subsec:parent}, we introduce the concept of slender genera and explain their connection to idoneal genera. The Smith--Minkowski--Siegel mass formula, a key tool for this section, is reviewed in §\ref{subsec:mass.formula} with notation drawn from Conway and Sloane's work~\cite{ConwaySloane88}. In §\ref{subsec:comparison.masses}, we develop a series of propositions comparing the masses of \(L \oplus [1]\) and \(L\), leading to \Cref{thm:bound.rank.>2}. 

The proof of \Cref{thm:idoneal} is presented in §\ref{subsec:finiteness}, and §\ref{subsec:enumeration} concludes this part with an algorithm to enumerate all idoneal genera of rank \(n \geq 3\).

\subsection{Slender genera} \label{subsec:parent}
Given a lattice \(L\), we denote by \(\Aut(L)\) its automorphism group, which is finite if \(L\) is positive definite.

Let \(g\) be a positive definite genus. Since each genus is a finite set (see, e.g., \cite[Kapitel VII, Satz (21.3)]{kneser}), one can define the \emph{mass} of \(g\) in the following way:
\[
    m(g) \coloneqq \sum_{i = 1}^h \frac{1}{|{\Aut(L^{(i)})}|},
\]
where \(L^{(1)},  L^{(2)}, \ldots, L^{(h)}\) is a complete set of representatives of lattices in \(g\).

\begin{definition} \label{def:parent}
    If \(L\) and \(M\) are two lattices in \(g\), then \(L \oplus [1]\) and \(M \oplus [1]\) belong to the same genus, which we denote by \(g'\). We say that \(g'\) is the \emph{child} of \(g\), and that \(g\) is a \emph{parent} of \(g'\).
\end{definition}

Consider the infinite graph in which each node represents a genus, and the nodes representing \(g\) and \(g'\) are joined by an edge. 
As explained in \cite[§10]{ConwaySloane88}, each connected component of this graph is called a \emph{vine} and consists of a single path, called the \emph{stem}, formed by the nodes representing odd genera, together with other nodes, called \emph{twigs}, representing even genera. Each twig is joined to the stem by a single edge, a fact which we restate in the following lemma.
\begin{lemma}[{cf. \cite[Lemma 3 on p. 280]{ConwaySloane88}}] \label{lem:at.most.two.parents}
    Each genus has at most two parents: possibly one odd parent on the stem, and possibly one even parent on a twig. 
\end{lemma} 

\begin{definition} \label{def:slender}
    We say that \(g\) is \emph{slender} if \(m(g') \leq m(g)\).
\end{definition}

The following easy, but crucial, lemma sets a connection between idoneal and slender genera.

\begin{lemma} \label{lem:slender.parent}
	If an idoneal genus \(g\) has exactly one parent \(f\), then
	\[
		2 m(g) \leq m(f).
	\]
	If an idoneal genus \(g\) has two parents \(f_1,f_2\), then
	\[
		2 m(g) \leq m(f_1) + m(f_2).
	\]	
	In particular, any idoneal genus has at least one slender parent.
\end{lemma}
\begin{proof}
    Let \(g\) be an idoneal genus. By definition, every lattice in \(g\) is of the form \(L \oplus [1]\) for some lattice \(L\) in a parent of \(g\). Therefore, \(g\) has at least one parent, and, by \Cref{lem:at.most.two.parents}, at most two.
    
    If \(g\) has only one parent \(f\), then all lattices in \(g\) are of the form \(L^{(i)} \oplus [1]\), where \(L^{(i)}\) runs over the lattices in \(f\). 
    By \cite[Satz~27.5]{kneser}, the positive definite lattices \(L^{(i)} \oplus [1]\) and \(L^{(j)} \oplus [1]\) are isomorphic if and only if \(i = j\). 
    Given that
    \begin{equation} \label{eq:Aut(L+1)}
    	|{\Aut(L \oplus [1])}| \geq 2 |{\Aut(L)}|,
    \end{equation}
    we obtain
    \[
        m(g) = \sum_{i = 1}^h \frac{1}{|{\Aut(L^{(i)} \oplus [1])}|} \leq \frac 12 \sum_{i = 1}^h \frac{1}{|{\Aut(L^{(i)})}|} = \frac 12 m(f),
    \]
    proving the first assertion. In particular, \(f\) must be a slender genus.
    
    Now, if \(g\) has two parents \(f_1,f_2\), then all lattices in \(g\) are either of the form \(L_1^{(i)} \oplus [1]\), where \(L_1^{(i)}\) runs over the lattices in~\(f_1\), or of the form \(L_2^{(j)} \oplus [1]\), where \(L_2^{(j)}\) runs over the lattices in~\(f_2\).
    Using \eqref{eq:Aut(L+1)} again, we have
    \begin{align*}
        m(g) = & \sum_{i = 1}^h \frac{1}{|{\Aut(L_1^{(i)} \oplus [1])}|} + \sum_{j = 1}^k \frac{1}{|{\Aut(L_2^{(j)} \oplus [1])}|} \\ 
        & \leq \frac 12 \sum_{i = 1}^h \frac{1}{|{\Aut(L_1^{(i)})}|} + \frac 12 \sum_{j = 1}^k \frac{1}{|{\Aut(L_2^{(j)})}|} = \frac 12 (m(f_1) + m(f_2)). 
    \end{align*} 
    If neither \(f_1\) nor \(f_2\) were slender, we would have \(m(f_1) < m(g)\) and \(m(f_2) < m(g)\), which would lead to a contradiction.
\end{proof}

\begin{example}
    If an idoneal genus has two parents, not both need to be slender.  
    Consider, for instance, the genus \(g\) of unimodular, positive definite lattices of rank~\(9\), which consists of two classes: \(g = \{9[1], \bE_8 \oplus [1]\}\), and has mass
    \[
        m(g) = \frac{1}{2^9 \cdot 9!} + \frac{1}{2 \cdot 696\,729\,600} = \frac{17}{2\,786\,918\,400} \approx 6.099 \cdot 10^{-9}.
    \]
    This genus is idoneal and has two parents: \(f_1 = \{8[1]\}\) and \(f_2 = \{\bE_8\}\), with masses
    \begin{align*}
        m(f_1) & = \frac{1}{2^8 \cdot 8!} = \frac{1}{10\,321\,920} \approx 9.688 \cdot 10^{-8}, \\  
        m(f_2) & = \frac{1}{696\,729\,600} \approx 1.435 \cdot 10^{-9}.
    \end{align*}
    By definition, \(f_1\) is slender, but \(f_2\) is not.
\end{example}



\subsection{The mass formula} \label{subsec:mass.formula}

We recall here a formula, called the Smith--Minkowski--Siegel mass formula, that allows us to compute the mass of a genus \(g\) starting from any lattice \(L\) in \(g\). 
We follow very closely the notation of Conway--Sloane \cite{ConwaySloane88} (but note that we work with lattices \(L\) instead of forms \(f\)).

Considering a prime $p$ and the $p$-adic lattice \(L\otimes \IZ_p\), we can fix a \(p\)-adic Jordan decomposition:
\begin{equation} \label{eq:p-adic.decomposition}
    L \otimes \IZ_p= \cdots \oplus L_{1/p}\big( \tfrac{1}{p} \big) \oplus L_1 \oplus L_p(p) \oplus L_{p^2}(p^2) \oplus \cdots = \bigoplus_q L_q(q),
\end{equation}
where \(q\) ranges over all powers of \(p\), including those with negative exponent,  each~\(L_q\) is a \(p\)-adically integral lattice whose determinant is prime to \(p\), and \(L_q(q)\) has the same underlying module as~\(L_q\), but the values of its quadratic form are multiplied by~\(q\).
All but finitely many lattices \(L_q\) have rank~\(0\), so the sum in \eqref{eq:p-adic.decomposition} is finite.

For \(p = 2\), a lattice \(L_q\) is called \emph{of type \(\I\)} or \emph{odd} if it represents an odd \(2\)-adic integer, otherwise it is called \emph{of type \(\II\)} or \emph{even}. Moreover, \(L_q\) is called \emph{bound} if either (or both) of the adjacent constituents \(L_{q/2}\) or \(L_{2q}\) is of type \(\I\), otherwise it is called \emph{free}. One also defines an invariant modulo \(8\), called the \emph{octane value} of \(L_q\), in the following way. 
If \(L_q\) has type \(\II\), then its octane value is \(0\) if \(\det(L_q) \equiv \pm 1 \bmod 8\), or \(4\) if \(\det(L_q) \equiv \pm 3 \bmod 8\). 
If \(L_q\) has type \(\I\), it is the orthogonal direct sum of \(2\)-adic lattices of rank $1$ with gram matrices given by \(2\)-adic units \(a_i\in \ZZ_2^\times\) (see, e.g., \cite[Lemma 4.1]{Allcock.Gal.Mark}). Then, its octane value is equal to the number of \(a_i\) that are congruent to \(1 \bmod 4\) minus the number of \(a_i\) that are congruent to \(-1 \bmod 4\).

Let \(n\) and \(d\) be the rank and determinant of~\(L\), respectively. Put \(n(q) \coloneqq \dim L_q\).
We henceforth assume that \(L\) is an integral lattice, so that \(n(q) = 0\) for \(q < 1\) and 
\begin{equation} \label{eq:n=sum_n(q)}
     n = \sum_q n(q).
\end{equation}
Computing determinants, we also infer that
\begin{equation} \label{eq:nu_p(d)}
    \sum_{q} \nu_p(q) n(q) = \nu_p(d).
\end{equation}
Combining \eqref{eq:n=sum_n(q)} with \eqref{eq:nu_p(d)}, we derive the following lower bound for the rank of \(L_1\):
\begin{equation} \label{eq:n(1)}
    \dim(L_1) = n(1) = n - \sum_{1 < q}n(q) \geq n - \sum_{1 \leq q} \nu_p(q) n(q) = n - \nu_p(d).
\end{equation}

The \emph{species} of \(L_q\) is a symbol which takes the form \(N\) (for \(N\) odd), \(0+\), \(N+\) or \(N-\) (for \(N\) even and positive). For \(p = 2\), the species of \(L_q\) depends on its type, rank, status (free or bound) and octane value. For odd \(p\), the species of \(L_q\) depends on \(p \bmod 4\), \(n \bmod 4\), and \(\det(L_q) \bmod p\). For the precise definition, see \cite[Table~1]{ConwaySloane88}.

Now, one defines the \emph{\(p\)-mass} of~\(L\), denoted by \(m_p(L)\), as the product of three factors:
\begin{equation} \label{eq:p-mass}
    m_p(L) \coloneqq \underbrace{\prod_q M_p(L_q)}_{\Delta_p(L)} \cdot \underbrace{\prod_{\substack{q, \, q' \\ q < q'}} (q'/q)^{\frac{1}{2} n(q)n(q')}}_{\chi_p(L)} \cdot \underbrace{2^{n(\I,\I) - n(\II)}}_{\substack{\tf_2(L) \\ (\text{only if } p = 2)}}
\end{equation}

The first factor, \(\Delta_p(L)\), is called \emph{diagonal product}, while the factor \(M_p(L_q)\) is called \emph{diagonal factor}. The value of \(M_p(L_q)\) depends on the species of \(L_q\) as follows:  
\begin{align*}
    M_p(0+) & \coloneqq 1, \\
    M_p(2s-1) & \coloneqq \frac{1}{2(1-p^{-2})(1-p^{-4})\cdots(1-p^{2-2s})} & (s > 0), \\
    M_p(2s\pm) & \coloneqq \frac{1}{2(1-p^{-2})(1-p^{-4})\cdots(1-p^{2-2s}) \cdot (1 \mp p^{-s})} & (s > 0).
\end{align*}

The second factor, \(\chi_p(L)\), is called \emph{cross-product}. Observe that 
\begin{equation} \label{eq:chi_p(L)=1}
    \chi_p(L) = 1 \quad \text{if \(p \nmid d\)},
\end{equation}
since, in this case, \(n(q) = 0\) for all \(q \neq 1\) due to \eqref{eq:nu_p(d)} and \(n(q)=0\) for \(q<1\) by the integrality of \(L\). 

The third and last factor \(\tf_2(L)\) is called the \emph{type factor} and is present only if \(p = 2\), in which case \(n(\II)\) is the sum of the ranks of all \(2\)-adic Jordan constituents~\(L_q\) of type \(\II\), and \(n(\I,\I)\) is the total number of pairs of adjacent constituents \(L_q,L_{2q}\) that are both of type \(\I\).

The following formula, which is valid for any lattice \(L\) in a positive definite genus~\(g\) of rank \(n \geq 2\), is known as the \emph{Smith--Minkowski--Siegel mass formula} \cite[p.~263, eq.~(2)]{ConwaySloane88}:
\begin{equation} \label{eq:smith.minkowski.siegel}
    m(g) = 2\pi^{-\frac{1}{4}n(n+1)} \prod_{j=1}^n \Gamma\big(\tfrac{1}{2}j\big) \prod_p 2m_p(L) \qquad (n \geq 2).
\end{equation}
We refer to \cite[(33.6), (35.1)]{kneser} for a modern proof and to \cite{cho2015} for a proof of the \(2\)-adic densities.  
The product in \eqref{eq:smith.minkowski.siegel} runs over all prime numbers \(p = 2, 3, 5, \ldots\).  
For \(n > 2\), it is absolutely convergent, but for \(n = 2\), this fails. However, if the primes are sorted in ascending order, the product converges (see, e.g., \cite[\S 109]{landau}).  
Thus, we will adopt this convention for all infinite products \(\prod_p\) over prime numbers.

\subsection{Comparison of masses} \label{subsec:comparison.masses}
Let \(g\) be a positive definite genus of rank \(n \geq 2\) and determinant~\(d\), and let \(L\) be any lattice in \(g\).
We abbreviate 
\[
    L' \coloneqq L \oplus [1],
\]
and we write \(g'\) for the child of \(g\), i.e., the genus of \(L'\). We wish to compare the mass of \(g'\) with the mass of~\(g\). 
It follows from the mass formula~\eqref{eq:smith.minkowski.siegel} that
\begin{equation} \label{eq:comparison.masses}
    \frac{m(g')}{m(g)} = \pi^{-\frac{1}{2}(n+1)} \Gamma\big(\tfrac{1}{2}(n+1)\big) \prod_p \frac{m_p(L')}{m_p(L)} \qquad (n \geq 2).
\end{equation}
Recall that the right-hand side is indeed (conditionally) convergent.
In order to estimate the last factor \(\prod_p m_p(L')/m_p(L)\), we fix a \(p\)-adic decomposition of~\(L\) as in \eqref{eq:p-adic.decomposition} for any prime \(p\). 
Then, \(L'\) has a \(p\)-adic decomposition \(L' = \bigoplus_q L'_q(q)\) with 
\begin{equation} \label{eq:L'_q}
    L'_1 = L_1 \oplus [1] \quad \text{and} \quad L'_q = L_q \text{ for \(q \neq 1\)}.
\end{equation}

By the definition of \(p\)-mass (see \eqref{eq:p-mass}), we can write
\begin{equation} \label{eq:ABCDE}
    \prod_p \frac{m_p(L')}{m_p(L)} 
    =       \underbrace{\frac{\Delta_2(L')}{\Delta_2(L)}}_A
    \cdot   \underbrace{\prod_{p \mid d, p\neq 2} \frac{\Delta_p(L')}{\Delta_p(L)}}_B
    \cdot   \underbrace{\prod_{p \nmid 2d} \frac{\Delta_p(L')}{\Delta_p(L)}}_C
    \cdot \underbrace{\prod_p \frac{\chi_p(L')}{\chi_p(L)}}_D
    \cdot \underbrace{\frac{\tf_2(L')}{\tf_2(L)}}_E.
\end{equation}

We now estimate the factors \(A\), \(B\), \(C\), \(D\), and \(E\). The product in factor \(D\) is finite due to \eqref{eq:chi_p(L)=1}. Notably, only the product in factor \(C\) is infinite, and converges because the left-hand side of the equality is already known to be convergent.

\begin{lemma} \label{lem:diag.factors}
    The following inequalities hold:
    \begin{align*}
        \frac 12 &\leq \frac{M_p(2t+1)}{M_p(2t+)} \leq 1,  &  \frac{M_p((2t+2)+)}{M_p(2t+)} & \geq \frac{1}{2}, &  \frac{M_p(2t+3)}{M_p(2t+1)} & \geq 1 & (t\geq 0); \\
         1 & \leq \frac{M_p(2t+1)}{M_p(2t-)} \leq 2, & \frac{M_p((2t+2)-)}{M_p(2t-)} & \geq 1, &
    \frac{M_p(2t-)}{M_p(2t+3)}&\geq 1 &  (t> 0);\\
    &\frac{M_p((2t+2)\pm)}{M_p(2t+1)} \geq \frac{1}{2},& \frac{M_p(2t+3)}{M_2(2t+)} &\geq 1 & &&(t \geq 0).
    \end{align*}
\end{lemma}
\begin{proof}
    We compute the ratios of the diagonal factors using their definition. For instance,
    \[
        \frac{M_p(2t+1)}{M_p(2t+)} = \frac{2(1-p^{-2})(1-p^{-4})\cdots(1-p^{2-2t})\cdot (1 - p^{-t})}{2(1-p^{-2})(1-p^{-4})\cdots(1-p^{2-2t})\cdot (1-p^{-2t})} = \frac{1}{1 + p^{-t}},
    \]
    and the inequalities follow easily.
\end{proof}

\begin{proposition}[Factor \(A\)] \label{prop:factor.A}
    Given a lattice \(L\), we have
    \[
        \frac{\Delta_2(L')}{\Delta_2(L)} \geq \begin{cases} \tfrac{1}{2} & \text{if \(L\) is odd,} \\ \tfrac{1}{8} & \text{if \(L\) is even.} \end{cases}        
    \]
\end{proposition}
\begin{proof}
    Let \(q\) be a power of \(2\) (possibly with negative exponent). Because of \eqref{eq:L'_q}, the rank, type (\(\I\) or \(\II\)) and octane value of \(L_q'\) for \(q \neq 1\) is the same as the rank, type and octane value of \(L_q\), respectively. 
    Moreover, the status (free or bound) of \(L_q'\) for \(q \neq \tfrac 12, 2\) is the same as the status of \(L_q\). Thus, \(L_q\) and \(L_q'\) have the same species, hence the same diagonal factor, for \(q \neq \tfrac 12, 1, 2\). It follows that
    \begin{equation} \label{eq:Delta2}
        \frac{\Delta_2(L')}{\Delta_2(L)} = \frac{\prod_q M_2(L_q')}{\prod_q M_2(L_q)} = 
        \frac{M_2(L_{1/2}')}{M_2(L_{1/2})} \cdot \frac{M_2(L_{1}')}{M_2(L_{1})} \cdot \frac{M_2(L_{2}')}{M_2(L_{2})}.
    \end{equation}
    
    Suppose first that \(L\) is odd, i.e., \(L_1\) is of type \(\I\). Then \(L_{1/2}\) and \(L_2\) are bound, and so are \(L_{1/2}'\), \(L_2'\). Therefore, \(M_2(L_{1/2}') = M_2(L_{1/2})\) and \(M_2(L_{2}') = M_2(L_{2})\). Since both the rank and the octane value of \(L_1'\) increase by \(1\) with respect to~\(L_1\), we infer the following by looking at \cite[Table~1]{ConwaySloane88}:
    \begin{itemize}
        \item if the species of \(L_1\) is \(2t+\), then the species of \(L_1'\) is either \(2t+\), \((2t+2)+\) \(2t+1\), or \(2t+3\);
        \item if the species of \(L_1\) is \(2t+1\), then the species of \(L_1'\) is either \(2t+\), \((2t+2)+\), \(2t+1\), \(2t+3\), \(2t-\), or \((2t+2)-\);
        \item if the species of \(L_1\) is \(2t-\), then the species of \(L_1'\) is either \(2t+1\), \(2t+3\), \(2t-\), or \((2t+2)-\).
    \end{itemize}
    In each case, it follows from \Cref{lem:diag.factors} that \(M_2(L_1')/M_2(L_1) \geq \frac 12 \). Hence, we have \(\Delta_2(L')/\Delta_2(L) \geq 1\cdot \frac 12 \cdot 1 = \frac 12\) from \eqref{eq:Delta2} as wished.

    Suppose now that \(L\) is even, i.e., \(L_1\) is of type \(\II\). Then, \(L_{1/2}'\) is bound, while \(L_{1/2}\) is free, so \(M_2(L_{1/2}') = \frac 12\) and \(M_2(L_{1/2}) = 1\). Also, \(L_2'\) is bound, so its species is \(2t+1\), while the species \(L_2\) can be either \(2t+\), \(2t+1\) or \(2t-\). In any case, \(M_2(L_2')/M_2(L_2) \geq \frac 12 \) by \Cref{lem:diag.factors}. Arguing as before with \(L_1'\) and \(L_1\), we infer again that \(M_2(L'_1)/M_2(L_1) \geq \frac 12\). All in all, equation \eqref{eq:Delta2} implies \(\Delta_2(L')/\Delta_2(L) \geq \frac 12 \cdot \frac 12 \cdot \frac 12 = \frac 18\), as desired.
\end{proof}

\begin{proposition}[Factor \(B\)] \label{prop:factor.B}
    For \(n, d \in \IZ_{\geq 1}\) and \(p\) a prime, define
    \begin{align*}
    	\mu(n,p,d) & \coloneqq \max\Big(0, \Big\lceil \frac{n-\nu_p(d)}{2}\Big\rceil\Big), \\
        \xi(n, d) & \coloneqq \prod_{p \mid d, p \neq 2} \frac{1}{1 + p^{-\mu(n,p,d)}}.
    \end{align*} 
    Given a positive definite lattice \(L\) of rank \(n\) and determinant~\(d\), we have
    \[
        \prod_{p \mid d, p \neq 2} \frac{\Delta_p(L')}{\Delta_p(L)} \geq \xi(n, d).
    \]
\end{proposition}
\begin{proof}
    Let \(p\) be an odd prime dividing \(d\). Because of \eqref{eq:L'_q}, \(L_q'\) has the same rank and determinant as \(L_q\) for \(q \neq 1\). Thus, \(L_q\) and \(L_q'\) have the same species, hence the same diagonal factor, for \(q \neq 1\). It follows that
    \[
        \frac{\Delta_p(L')}{\Delta_p(L)} = \frac{\prod_q M_p(L_q')}{\prod_q M_p(L_q)} = \frac{M_p(L_1')}{M_p(L_1)}.
    \]
    
    Suppose first that \(\dim L_1 = 2s\) with \(s \geq 0\). By looking at \cite[Table~1]{ConwaySloane88}, we see that the species of \(L_1\) is either \(2s+\) or \(2s-\), while the species of \(L_1'\) is \(2s+1\). Therefore, we have (also for \(s = 0\))
    \[
        \frac{M_p(L_1')}{M_p(L_1)} = \frac{M_p(2s+1)}{M_p(2s\pm)} = \frac{1 \mp p^{-s}}{1 - p^{-2s}} = \frac{1}{1 \pm p^{-s}} \geq \frac{1}{1 + p^{-s}}.
    \]

    Suppose now that \(\dim L_1 = 2s-1\) with \(s > 0\). Looking at \cite[Table~1]{ConwaySloane88}, we see that the species of \(L_1\) is \(2s-1\), while the species of \(L_1'\) is either \(2s+\) or \(2s-\). In both cases, we come to the same bound as above:
    \[
        \frac{M_p(L_1')}{M_p(L_1)} = \frac{M_p(2s\pm)}{M_p(2s+1)} = \frac{1}{1 \mp p^{-s}} \geq \frac{1}{1 + p^{-s}}.
    \]
    
    We observe that inequality \eqref{eq:n(1)} implies
    \[
        s \geq \mu(n,p,d).
    \]
    The result follows by multiplying over all odd prime divisors of \(d\).
\end{proof}

Following \cite[§7]{ConwaySloane88}, we introduce the function \(\zeta_D\) for \(s > 1\) and \(D \in \mathbb{Z}\), or \(s=1\) and \(D \in \ZZ_{<0}\):
\[
    \zeta_D(s) \coloneqq \prod_p \frac{1}{1 - (D/p)p^{-s}} = \sum_{2\nmid n} (D/n)n^{-s},
\]
where \((D/p)\) is a Legendre symbol, to be interpreted as \(0\) if \(p \mid 2D\), and \((D/n)\) is a Jacobi symbol. 
The product runs over all prime numbers in ascending order. 
We remark that 
both sides converge absolutely for \(s > 1\) and conditionally for \(s = 1\) and $D<0$. We refer to \cite[\S 109, p.~446, p.~449]{landau} and \cite[p.~153]{kneser}) for the equality and conditional convergence (see also \cite[§12.7]{hua}). 


\begin{proposition}[Factor \(C\)] \label{prop:factor.C}
Let \(L\) be a positive definite lattice of rank \(n \geq 2\) and determinant~\(d\). Write \(n = 2s-1\) or \(n = 2s\) with \(s \in \IZ_{\geq 1}\). 
\begin{enumerate}
    \item[(a)] If \(s > 1\), then
\[
    \prod_{p \nmid 2d} \frac{\Delta_p(L')}{\Delta_p(L)} \geq (1 + 2^{-s}) \frac{\zeta(2s)}{\zeta(s)}. 
\]

\item[(b)] If \(s = 1\), then 
\[\prod_{p \nmid 2d} \frac{\Delta_p(L')}{\Delta_p(L)} \geq \frac{1}{\zeta_{-d}(1)}.\]
\end{enumerate}
\end{proposition}
\begin{proof}
    If \(p \nmid 2d\), then the \(p\)-adic Jordan decomposition of~\(L\) is simply \(L = L_1\). Therefore, \(m_p(L) = \Delta_p(L) = M_p(L_1)\) takes on the so-called \emph{standard value} 
    \[
        \std_p(L) \coloneqq \frac{1}{2(1-p^{-2})(1-p^{-4})\cdots (1-p^{2-2s})\cdot (1-\epsilon p^{-s})},
    \]
    where \(\epsilon\) is \(0\) for \(n\) odd, and is otherwise the Legendre symbol \((D/p)\) with \(D = (-1)^s d\).

    If \(L\) is of rank \(n = 2s\), then \(L'\) is of rank \(2s+1\), hence
    \[
        \frac{\Delta_p(L')}{\Delta_p(L)} = \frac{\std_p(L')}{\std_p(L)} = \frac{1 - (D/p) p^{-s}}{1 - p^{-2s}}=\frac{1}{1+(D/p)p^{-s}}.
    \]
    If \(L\) is of rank \(n = 2s-1\), then \(L'\) is of rank \(2s\), hence
    \[
        \frac{\Delta_p(L')}{\Delta_p(L)} = \frac{\std_p(L')}{\std_p(L)} = \frac{1}{1 - (D/p) p^{-s}}.
    \]
    Multiplying over all primes \(p \nmid 2d\) we have
    \[
        \prod_{p \nmid 2d} \frac{\Delta_p(L')}{\Delta_p(L)} = \prod_{p \nmid 2d} \frac{1}{1 \pm (D/p) p^{-s}}.
    \]

    Now suppose \( s > 1 \). Using the inequality \( 1 \pm (D/p)p^{-s} \leq 1 + p^{-s} \), we can estimate as follows:
    \begin{equation*}
    \prod_{p \nmid 2d} \frac{\Delta_p(L')}{\Delta_p(L)} 
    = \prod_{p \nmid 2d} \frac{1}{1 \pm (D/p)p^{-s}} 
    \geq \prod_{p \nmid 2d} \frac{1}{1 + p^{-s}} 
    \geq \prod_{p \nmid 2} \frac{1}{1 + p^{-s}}.
    \end{equation*}
    This simplifies further:
    \begin{equation*}
    \prod_{p \nmid 2} \frac{1}{1 + p^{-s}} 
    = (1 + 2^{-s}) \prod_{p} \frac{1}{1 + p^{-s}} \\
    = (1 + 2^{-s}) \frac{\zeta(2s)}{\zeta(s)},
\end{equation*}
where the condition \( s > 1 \) ensures the convergence of all the infinite products involved.

Next, consider the case \( s = 1 \). In particular, we have \(n = 2\) and \(D = -d\). Observing that \( 1 - p^{-2} \leq 1 \) and noting that \( (D/p) = 0 \) for \( p \mid 2d \), we obtain:
\begin{align*}
    \prod_{p \nmid 2d} \frac{\Delta_p(L')}{\Delta_p(L)} 
    = \prod_{p \nmid 2d} \frac{1 - (D/p)p^{-1}}{1 - p^{-2}} 
    \geq \prod_{p \nmid 2d} (1 - (D/p)p^{-1}) 
    = \frac{1}{\zeta_D(1)}.
\end{align*}
Here, recall that the products in the argument are conditionally convergent. Absolute convergence is not required because the order of the factors is not altered.

\end{proof}

\begin{proposition}[Factor \(D\)] \label{prop:factor.D}
    Given a positive definite lattice \(L\) of determinant~\(d\), we have
    \begin{equation*} 
        \prod_p \frac{\chi_p(L')}{\chi_p(L)} = \sqrt{d}.
    \end{equation*}
\end{proposition}
\begin{proof}
    Denoting \(n'(q) \coloneqq \dim L'_q\) for any power \(q\) of \(p\), we have \(n'(1) = n(1) + 1\) and \(n'(q) = n(q)\) for \(q \neq 1\) as a consequence of \eqref{eq:L'_q}. Recall that \(n(q) = 0\) for \(q<1\) by the integrality of~\(L\). By definition, the cross-product of~\(L'\) is given by
    \begin{align*}
        \chi_p(L') & = \prod_{1 = q < q'} (q'/q)^{\frac{1}{2} n'(q) n'(q')} \prod_{1 < q < q'} (q'/q)^{\frac{1}{2} n'(q) n'(q')} \\
                & = \prod_{1 < q'} (q')^{\frac 12 (n(1)+1)n(q')} \prod_{1 < q < q'} (q'/q)^{\frac{1}{2} n(q) n(q')} \\
                & = \prod_{1 < q'} (q')^{\frac 12 n(q')} \prod_{1 \leq q < q'} (q'/q)^{\frac{1}{2} n(q) n(q')} \\
                & = p^{\frac{1}{2}\nu_p(d)} \chi_p(L),
    \end{align*}
    where the last equality follows from \eqref{eq:nu_p(d)}. The result is obtained by multiplying over all primes.
\end{proof}

\begin{proposition}[Factor \(E\)] \label{prop:factor.E}
For a lattice \(L\) of rank \(n\) and determinant~\(d\), the following inequality holds:
\[
    \frac{\tf_2(L')}{\tf_2(L)} \geq 
        \begin{cases} 
            1 & \text{if \(L\) is odd,} \\
            2^{\max(0,n-\nu_2(d))} & \text{if \(L\) is even.}
        \end{cases}
\]
\end{proposition}
\begin{proof}
    Fix a \(2\)-adic decomposition of~\(L\) as in \eqref{eq:p-adic.decomposition} and consider the induced \(2\)-adic decomposition of \(L'\) as in \eqref{eq:L'_q}. Let \(n'(\II)\) be the sum of the ranks of all \(2\)-adic Jordan constituents \(L_q'\) that have type \(\II\), and \(n'(\I,\I)\) be the total number of pairs of adjacent constituents \(L_q',L_{2q}'\) that are both of type \(\I\). 
    
    If \(L\) is odd, then \(L_1\) is of type \(\I\) and so is \(L_1'\). 
    Hence, we have \(n'(\II) = n(\II)\) and \(n'(\I,\I) = n(\I,\I)\), i.e., \(\tf_2(L') = \tf_2(L)\).

    If \(L\) is even, then \(L_1\) is of type \(\II\), but \(L_1'\) is of type \(\I\). Hence, we have \(n'(\II) = n(\II) - \dim L_1\) and \(n'(\I,\I) \geq n(\I,\I)\), so
    \[
        \tf_2(L') = 2^{n'(\I,\I) - n'(\II)} \geq 2^{n(\I,\I) - n(\II) + \dim L_1} = 2^{\dim L_1} \tf_2(L).
    \]
    We conclude by observing that \(\dim L_1 \geq \max(0,n-\nu_2(d))\) because of \eqref{eq:n(1)}.
\end{proof}

\begin{theorem} \label{thm:bound.rank.>2}
For \(n, d \in \IZ_{\geq 1}\), and $\xi(n,d)$ as in \Cref{prop:factor.B}, define the following functions:
\[
\Phi_n^\I(d) \coloneqq \tfrac{1}{2} \xi(n,d) \sqrt{d}, \quad  
\Phi_n^\II(d) \coloneqq \tfrac{1}{8} \xi(n,d) \sqrt{d} \cdot 2^{\max(0, n - \nu_2(d))}.
\]
Assume \(n \geq 3\) and write \(n = 2s - 1\) or \(n = 2s\), where \(s \in \IZ_{\geq 2}\). Define
\[
c_n \coloneqq \pi^{-\frac{1}{2}(n+1)} \Gamma\big(\tfrac{1}{2}(n+1)\big) (1 + 2^{-s}) \frac{\zeta(2s)}{\zeta(s)}.
\]
For a positive definite genus \(g\) of rank \(n\) and determinant \(d\), the following inequality holds:
\[
\frac{m(g')}{m(g)} \geq  
\begin{cases}  
    c_n \Phi_n^\I(d) & \text{if \(L\) is odd,} \\  
    c_n \Phi_n^\II(d) & \text{if \(L\) is even.}  
\end{cases}
\]
\end{theorem}
\begin{proof}
	We express \(m(g')/m(g)\) as in \eqref{eq:comparison.masses}. 
    In \Cref{prop:factor.A,prop:factor.B,prop:factor.C,prop:factor.D,prop:factor.E}, we derived lower bounds for the factors \(A, B, C, D, E\) appearing in \eqref{eq:ABCDE}. The functions \(\Phi_n^I\) and \(\Phi_n^{II}\) encapsulate the estimates for \(A, B, D,\) and \(E\), while the estimate for \(C\) is absorbed into \(c_n\). 
    The statement of the theorem then follows directly.  
    Finally, note that the hypothesis \(n \geq 3\) is required to apply \Cref{prop:factor.C}(a).
\end{proof}

From the definition of slender genus (\Cref{def:slender}), we infer the following corollary.

\begin{corollary} \label{cor:DIn}
Let \(n \geq 3\). With the notation of \Cref{thm:bound.rank.>2}, define
\begin{align*}
    D_n^\I & \coloneqq \{ d \in \IZ_{\geq 1} \mid c_n \Phi_n^\I(d) \leq 1 \}, \\
    D_n^\II & \coloneqq \{ d \in \IZ_{\geq 1} \mid c_n \Phi_n^\II(d) \leq 1 \}.
\end{align*}
If \(g\) is an odd (resp. even) slender genus of rank \(n\) and determinant \(d\), then \(d \in D^\I_n\) (resp. \(d \in D^\II_n\)).
\end{corollary}

Next, we turn to the case \(n = 2\) (equivalently, \(s = 1\)).

\begin{lemma} \label{lem:bound.zeta_D}
For all \(d > 0\), we have
\[
    \zeta_{-d}(1) \leq \tfrac{3}{2} \ln(4d) + 3.
\]
\end{lemma}
\proof
Put \(D = -d\) and let \(e \equiv 0,1 \bmod 4\) where \(e\) is not a perfect square. 
Define \(\psi_e(m)\) as the Kronecker symbol~\(\left(\frac{e}{m} \right)\). The associated \(L\)-series is 
 \(L(\psi_e,s) = \sum_{m=1}^\infty \psi_e(m)m^{-s}\). 
By \cite[§12.14, Theorem~14.3]{hua} we have \(L(\psi_e,1) \leq 2 + \ln|e|\).

Now, if \(D\equiv 0 \bmod 4\), then \(\zeta_D(s) = L(\psi_D,s)\).
If \(D\equiv 2,3 \bmod 4\), then \(\zeta_D(s)=\zeta_{4D}(s) = L(\psi_{4D},s)\).
Finally, if \(D\equiv 1 \bmod 4\), then 
\begin{align*}
(1-\psi_D(2)2^{-s})L(\psi_{D},s)
    &= \sum_{m=1}^\infty \psi_{D}(m)m^{-s} - \sum_{m=1}^\infty \psi_{D}(2m)(2m)^{-s}\\
    &= \sum_{m=1,3,5,\ldots}^\infty \psi_{D}(m)m^{-s}\\
    &= \zeta_D(s).
\end{align*}
Note that the series involved here converge conditionally for $s\geq 1$ because $L(\psi_{D},s)$ does.
Absolute convergence is not used in the argument, because we do not change the summation order.
Note that \((1-\psi_D(2)2^{-1})\leq 3/2\).  
In any case, we conclude that
\[
\zeta_D(1) \leq \tfrac{3}{2} \ln(4d) + 3. \qedhere
\]
\endproof 

\begin{theorem} \label{thm:bound.rank.2}
Define \(\tilde c_2 \coloneqq \pi^{-\frac{3}{2}} \Gamma\left(\tfrac{3}{2}\right) = \tfrac{1}{2\pi}\). For \(d \geq 1\), let  
\[
\tilde{\Phi}_2^\I(d) \coloneqq \frac{\Phi_2^\I(d)}{\tfrac{3}{2} \ln(4d) + 3}, \qquad  
\tilde{\Phi}_2^\II(d) \coloneqq \frac{\Phi_2^\II(d)}{\tfrac{3}{2} \ln(4d) + 3},
\]  
where \(\Phi_2^\I\) and \(\Phi_2^\II\) are as defined in \Cref{thm:bound.rank.>2}.  
For a positive definite genus \(g\) of rank \(n = 2\) and determinant \(d\), the following inequalities hold:  
\[
\frac{m(g')}{m(g)} \geq  
\begin{cases}  
\tilde c_2 \Phi_2^\I(d)/\zeta_{-d}(1) \geq \tilde c_2 \tilde{\Phi}_2^\I(d) & \text{if \(g\) is odd,} \\  
\tilde c_2 \Phi_2^\II(d)/\zeta_{-d}(1) \geq \tilde c_2 \tilde{\Phi}_2^\II(d) & \text{if \(g\) is even.}  
\end{cases}
\]
\end{theorem}
\proof
The statement follows from \eqref{eq:comparison.masses} and \eqref{eq:ABCDE}, by comparing the factors \(A,\ldots,E\) using respectively \Cref{prop:factor.A,prop:factor.B}, 
\Cref{prop:factor.C}(b), \Cref{prop:factor.D,prop:factor.E}. 
The functions \(\Phi^I_2\) and \(\Phi^{II}_2\) correspond to \(ABDE\) whereas $C\geq 1/\zeta_{-d}(1)$. The factor \(\tilde c_2\) is just the factor \(\pi^{-\frac{1}{2}(n+1)}\Gamma(\frac{1}{2}(n+1))\) with \(n = 2\).
The last inequalities follow from \Cref{lem:bound.zeta_D}.
\endproof

\begin{corollary} \label{cor:DI2}
With the notation of \Cref{thm:bound.rank.2}, define the following sets: 
\begin{align*}
    D_2^\I & \coloneqq \{d \in \IZ_{\geq 1} \mid \tilde{c}_2\Phi_2^\I(d)/\zeta_{-d}(1)\leq 1\},\\
    D_2^\II & \coloneqq \{d \in \IZ_{\geq 1} \mid \tilde{c}_2 \Phi_2^\II(d)/\zeta_{-d}(1)\leq 1\}, \\
    \tilde D_2^\I & \coloneqq \{d \in \IZ_{\geq 1} \mid \tilde{c}_2 \tilde{\Phi}_2^\I(d) \leq 1\}, \\
    \tilde D_2^\II & \coloneqq \{d \in \IZ_{\geq 1} \mid \tilde{c}_2 \tilde{\Phi}_2^\II(d) \leq 1\}.
\end{align*}
The following inclusions hold: \(D_2^\I \subseteq \tilde{D}_2^\I\) and \(D_2^\II \subseteq \tilde{D}_2^\II\).
Moreover, if \(g\) is an odd (resp. even) slender genus of rank \(2\) and determinant~\(d\), then \(d \in D^\I_2\) (resp. \(d \in D^\II_2\)).
\end{corollary}

\subsection{Finiteness} \label{subsec:finiteness}
Our goal is to establish the finiteness of slender genera of rank \(n \geq 2\), which will subsequently allow us to deduce the finiteness of idoneal genera.

\begin{lemma} \label{lem:Phi_new}
Let \(\Phi \colon \mathbb{Z}_{\geq 1} \to \mathbb{R}_{>0}\) and \(b \in \IR\). Define
\[
    D \coloneqq \{d \in \mathbb{Z}_{\geq 1} \mid \Phi(d) \leq b\}.
\]
Then, \(D\) is finite if there exist functions $m$ and $\varphi$ defined on the set of prime numbers and a prime $p_0$ such that the following conditions hold: 
\begin{enumerate}[(a)]
    \item \label{cond:Phi_(a)} For all prime numbers \(p\): \(m_p \in \IZ_{\geq 0}\), $\varphi(p)>1$ and
    \[
        \frac{\Phi(pd)}{\Phi(d)} \geq \varphi(p) \quad \text{for all \(d \in \mathbb{Z}_{\geq 1}\) with \(\nu_p(d) \geq m_p\)}.
    \]
    \item \label{cond:Phi_(b)} \(m_p = 0\) for all \(p \geq p_0\).
    \item \label{cond:Phi_(c)} The function \(\varphi(p)\) is increasing for \(p \geq p_0\) and satisfies
    \[
        \lim_{p \to \infty} \varphi(p) = +\infty.
    \]
\end{enumerate}
If moreover $m$, $\varphi$ and $p_0$ are computable, there exists an algorithm to compute \(D\).
\end{lemma}

\begin{proof}
    Define
    \[
        D_0 \coloneqq \{d \in \IZ_{\geq 1} \mid \text{\(\nu_p(d) \leq m_p\) for all primes \(p\)}\}.
    \]
    By condition \ref{cond:Phi_(b)}, \(D_0\) is finite and non-empty since \(1 \in D_0\).
    From condition \ref{cond:Phi_(a)}, 
    \begin{equation*}
        \Phi(pd) > \Phi(d)
    \end{equation*}
    for \(d\) with \(\nu_p(d) \geq m_p\), ensuring that \(\Phi(d)\) attains its minimum, \(\Phi_0\), within \(D_0\). In particular, if \(D \cap D_0 = \emptyset\), then \(D = \emptyset\). Assume now that \(D \neq \emptyset\).

    By condition \ref{cond:Phi_(c)}, there exists a prime \(p_1 \geq p_0\) such that
    \[
        \varphi(p_1) \Phi_0 > b.
    \]
    Thus, no prime \(p \geq p_1\) can divide any \(d \in D\). To see this, suppose \(d = pd'\) with \(p \geq p_1\). Then
    \[
        b \geq \Phi(pd') \geq \varphi(p) \Phi(d') \geq \varphi(p) \Phi_0 > b,
    \]
    a contradiction. Therefore, the primes dividing \(p \in D\) are bounded by \(p_1\).

    Let \(c \coloneqq \min_{p} \varphi(p)\).
    From condition \ref{cond:Phi_(a)} and \ref{cond:Phi_(c)}, we have \(c > 1\). For any \(p < p_1\), suppose \(p^m d' \in D\) for \(m \geq m_p\) and \(d' \in \IZ_{\geq 1}\).  By condition \ref{cond:Phi_(a)}, 
    \[
        b \geq \Phi(p^m d') \geq \varphi(p)^{m-m_p}\Phi(p^{m_p}d') \geq c^{m - m_p} \Phi_0. 
    \]  
    Rearranging gives  
    \[
        m \leq \frac{\ln(b) - \ln(\Phi_0)}{\ln(c)} + m_p.
    \]  

    Thus, if $d \in D$, the exponents \(\nu_p(d)\) for all \(p\) are bounded. Since both the primes dividing~\(d\) and their exponents are bounded, \(D\) is finite.

    Finally, let \(D'\) be the set produced by the following algorithm, where the function \verb|next_prime| applied to a positive integer \(a\) returns the smallest prime larger than~\(a\).

    \begin{myalgorithm}\caption{Determinant bounds}\label{alg:compute_D}
    \begin{algorithmic}[1]
        \REQUIRE \(\Phi\), \(b\), \(m_p\), \(\varphi(p)\), \(p_0\) as in the statement. \\
        \ENSURE \(D\) as in the statement.
    	\STATE \label{line:D0} \(D \longleftarrow \{d \mid \nu_p(d) \leq m_p \text{ for all primes \(p\)}, \Phi(d) \leq b\}\)
        \IF{\(D = \emptyset\)}
    		\RETURN \(\emptyset\)
    	\ENDIF
        \STATE \(\Phi_0 \longleftarrow \min_{d \in D} \Phi(d) \)
        \STATE \(p \longleftarrow 2\)
        \WHILE{\(p < p_0\) or \(\varphi(p)\Phi_0 \leq b\)} \label{state:while_p}
    	    \STATE \(E \longleftarrow D\)
            \WHILE{\(E \neq \emptyset\)} \label{state:while_nu}
                \STATE \(E \longleftarrow \{pe \mid e \in E, \Phi(pe) \leq b\}\)
                \STATE \(D \longleftarrow D \cup E\)
           \ENDWHILE
           \STATE \(p \longleftarrow \mathtt{next\_prime}(p)\)
        \ENDWHILE
        \RETURN \(D\)
    \end{algorithmic}
    \end{myalgorithm}
    
    We claim that \(D = D'\). Trivially, \(D' \subseteq D\), since each time that some \(d\) is added to \(D'\), the condition \(\Phi(d) \leq b\) has been checked. 
    
    In line \ref{line:D0}, the set \(D \cap D_0\) is computed. 

    If a prime \(p \geq p_0\) divides some \(d \in D\), then \(\varphi(p)\Phi_0 \leq b\), which justifies the condition on line \ref{state:while_p}. 
    
    Let \(d \in D\) with \(\nu_p(d) \geq m_p\). By induction on \(\sum_p (\nu_p(d) - m_p)\)
    (line~\ref{state:while_nu}) and on the smallest prime \(p\) appearing in the decomposition of \(d\) with \(\nu_p(d) > m_p\) (line~\ref{state:while_p}), we have \(d \in D'\). Thus, \(D \subseteq D'\), concluding the proof.
\end{proof}

In the following proposition, we use the notations of \Cref{thm:bound.rank.>2,thm:bound.rank.2} and \Cref{cor:DIn,cor:DI2}.

\begin{proposition} \label{prop:long}
    Let \(n\geq 2\). The following hold:
    \begin{itemize}
        \item The function \(\Phi(d) = \Phi_n^\I(d)\) satisfies the assumptions of \Cref{lem:Phi_new} with \(m_p = 0\) for all primes \(p\), \(p_0 = 3\) and
        \[
            \varphi(2) = \sqrt{2}, \quad \varphi(p) = \frac{p+1}{2p} \sqrt{p} \quad (p \geq 3).
        \]
        In particular, \(\Phi_n^\I(d)\) achieves its minimum at \(d = 1\). 
        \item The function \(\Phi(d) = \Phi_n^\II(d)\) satisfies the assumptions of \Cref{lem:Phi_new} with \(m_2 = n\), \(m_p = 0\) for all primes \(p \neq 2\), \(p_0 = 3\) and the same values of \(\varphi(p)\) as \(\Phi_n^\I(d)\).
        Moreover, \(\Phi_n^\II(d)\) achieves its minimum at \(d = 2^n\). 
        \item The function \(\Phi(d) = \tilde{\Phi}_2^\I\) satisfies the assumptions of \Cref{lem:Phi_new} with \(m_3 = 2\), \(m_p = 0\) for all primes \(p \neq 3\), \(p_0 = 11\) and
        \begin{align*}
            \varphi(2) &= \sqrt{2} \cdot \frac{\ln(4)+2}{\ln(8)+2}, & \varphi(3) &= \sqrt{3} \cdot \frac{\ln(36)+2}{\ln(108)+2}, \\ \varphi(5) & = \frac{3}{5} \sqrt{5} \cdot \frac{\ln(20)+2}{\ln(100)+2}, & 
            \varphi(7) & = \frac{4}{7} \sqrt{7} \cdot \frac{\ln(28)+2}{\ln(196)+2},\\
            \varphi(p) &= \frac{p+1}{2p} \sqrt{p} \cdot \frac{\ln(4)+2}{\ln(4p)+2} & (p \geq 11).
        \end{align*}
        \item The function \(\Phi(d) = \tilde{\Phi}_2^\II\) satisfies the assumptions of \Cref{lem:Phi_new} with \(m_2 = 2\), \(m_3 = 2\) and \(m_p = 0\) for all primes \(p \geq 5\), \(p_0 = 11\) and the same values of \(\varphi(p)\) as \(\tilde{\Phi}_2^\I(d)\).
    \end{itemize}
    For all \(n \geq 2\), the sets \(D_n^\I\) and \(D_n^\II\), as well as the sets \(\tilde D_2^\I\) and \(\tilde D_2^\II\), are finite.
\end{proposition}

\begin{proof}
    As a preliminary observation, we note that
    \begin{equation} \label{eq:1>=x/x+1>=x+1/2x}
        1 \geq \frac{x}{x+1} \geq \frac{x+1}{2x} \qquad (x \geq 3).
    \end{equation}
    We begin by examining \(\Phi(d) = \Phi^\I_n(d)\). By definition,
    \begin{equation*}
        \frac{\Phi_n^\I(pd)}{\Phi_n^\I(d)} = \frac{\xi(n,pd)}{\xi(n,d)} \sqrt{p}.
    \end{equation*}
    For \(p = 2\), \(\xi(n,2d) = \xi(n,d)\) implies
    \begin{equation} \label{eq:PhiI.p=2}
        \frac{\Phi_n^\I(2d)}{\Phi_n^\I(d)} = \sqrt{2} > 1.
    \end{equation}
    Now suppose \(p \geq 3\). First, assume \(\nu_p(d) = 0\), so that \(\nu_p(pd) = 1\). Then, 
    \begin{equation*}
        \frac{\Phi_n^\I(pd)}{\Phi_n^\I(d)} = \frac{1}{1+p^{-\mu(n,p,pd)}}\sqrt{p}.
    \end{equation*}
    Since \(n \geq 2\) by hypothesis, we find 
    \[
        \mu(n,p,pd) = \max\Big(0, \Big\lceil \frac{n-\nu_p(pd)}{2} \Big\rceil\Big) = \Big\lceil \frac{n-1}{2} \Big\rceil \geq 1,
    \]
    implying
    \begin{equation} \label{eq:PhiI.p>=3,nup(d)=0}
        \frac{\Phi_n^\I(pd)}{\Phi_n^\I(d)} \geq \frac{1}{1+p^{-1}} \sqrt{p}  = \frac{p}{p+1} \sqrt{p}. \qquad (p \geq 3, \, \nu_p(d) = 0).
    \end{equation}
    If \(\nu_p(d) > 0\), then
    \begin{equation} \label{eq:PhiI.p>=3,nu>0(a)}
        \frac{\Phi_n^\I(pd)}{\Phi^I_n(d)} = \frac{1+p^{-\mu(n,p,d)}}{1+p^{-\mu(n,p,pd)}} \sqrt{p}, \qquad (p \geq 3, \, \nu_p(d) > 0).
    \end{equation}
    Here, \(\mu(n,p,pd)\) can equal \(\mu(n,p,d)\) or \(\mu(n,p,d)-1\). In the first case, we have
    \[
        \frac{\Phi_n^\I(pd)}{\Phi_n^\I(d)} = \sqrt{p}. 
    \]
    In the second case, where necessarily \(\mu(n,p,d) \geq 1\), we have 
    \[
        \frac{\Phi_n^\I(pd)}{\Phi_n^\I(d)} = \frac{1+p^{-\mu(n,p,d)}}{1+p^{-\mu(n,p,d)+1}} \sqrt{p} \geq \frac{1+p^{-1}}{1+p^{-1+1}} \sqrt{p} = \frac{p+1}{2p} \sqrt{p} 
    \]
    because, for fixed \(p\), the function
    \(
        \frac{1+p^{-x}}{1+p^{-x+1}}
    \)
    is increasing for all \(x\). Using \eqref{eq:1>=x/x+1>=x+1/2x}, we find independently from the value of \(\mu(n,p,pd)\):
    \begin{equation} \label{eq:PhiI.p>=3,nup(d)>0}
        \frac{\Phi_n^\I(pd)}{\Phi_n^\I(d)} \geq \frac{p+1}{2p} \sqrt{p} \qquad (p \geq 3, \nu_p(d) > 0).
    \end{equation}
    From \eqref{eq:1>=x/x+1>=x+1/2x}, \eqref{eq:PhiI.p>=3,nup(d)=0} and \eqref{eq:PhiI.p>=3,nup(d)>0}, it follows
    \begin{equation*} \label{eq:PhiI.p>=3}
        \frac{\Phi_n^\I(pd)}{\Phi_n^\I(d)} \geq \frac{p+1}{2p} \sqrt{p} \qquad (p \geq 3).
    \end{equation*}
    Note that the function \(\frac{x+1}{2x}\sqrt{x}\) is increasing for \(x>1\), so 
    \[
        \varphi(p) = \frac{p+1}{2p} \sqrt{p} \geq \frac{3+1}{2\cdot 3} \sqrt{3} > 1 \qquad (p\geq 3).
    \]
    From \Cref{lem:Phi_new} and its proof, it follows that \(D_n^\I\) is finite for all \(n \geq 3\), and that \(\Phi_n^\I(d)\) achieves its minimum for \(d = 1\).

    We now consider \(\Phi(d)= \Phi_n^\II(d)\). By definition of \(\Phi_n^\II\), we have  
    \[
        \frac{\Phi_n^\II(pd)}{\Phi_n^\II(d)} = \frac{\Phi_n^\I(pd)}{\Phi_n^\I(d)} \cdot \frac{2^{\max(0, n - \nu_2(pd))}}{2^{\max(0, n - \nu_2(d))}}.
    \]  
    If \(\nu_2(d) \geq n\), then \(\nu_2(2d) \geq n+1\), so we have  
    \[
    \frac{\Phi_n^\II(pd)}{\Phi_n^\II(d)} = \frac{\Phi_n^\I(pd)}{\Phi_n^\I(d)} \qquad (\nu_2(d)\geq n \mbox{ or } p\geq 3),
    \]  
    and the same inequalities as above apply. By \Cref{lem:Phi_new}, the set \(D_n^\II\) is finite for all \(n \geq 3\).
    For \(\nu_2(d) < n\) we also have, by \eqref{eq:PhiI.p=2}, that
    \[
    \frac{\Phi_n^\II(2d)}{\Phi_n^\II(d)} = 
    \frac{1}{\sqrt{2}} <1 \qquad (\nu_2(d) <n).
    \] 
    The minimum value of \(\Phi_n^\II(d)\) is therefore achieved for \(d_{\min} = 2^n\).

    To prove the finiteness of the sets \(D_2^\I\) and \(D_2^\II\), it suffices to prove the finiteness of \(\tilde D_2^\I\) and \(\tilde D_2^\II\). Therefore, we now turn to \(\Phi(d) = \tilde{\Phi}_2^\I(d)\). By definition, we have
    \[
        \frac{\tilde{\Phi}_2^\I(pd)}{\tilde{\Phi}_2^\I(d)} 
        = \frac{\Phi_2^\I(pd)}{\Phi_2^\I(d)} \cdot \frac{\ln(4d)+2}{\ln(4pd)+2}.
    \]
    Computing the derivative, we see that, for fixed \(p\), the function
    \[
        x \longmapsto \frac{\ln(4x)+2}{\ln(4px)+2}
    \]
    is increasing for \(x > 0\). By \eqref{eq:PhiI.p=2}, since \(d \geq 1\), we have
    \[
        \frac{\tilde{\Phi}_2^\I(2d)}{\tilde{\Phi}_2^\I(d)} \geq \sqrt{2} \cdot \frac{\ln(4)+2}{\ln(4\cdot 2)+2} \approx 1.173 >1.
    \]
    Assume \(p \geq 3\). By \eqref{eq:PhiI.p>=3,nup(d)=0},
    \begin{equation} \label{eq:Phi2.nu=0}
        \frac{\tilde{\Phi}_2^\I(pd)}{\tilde{\Phi}_2^\I(d)} \geq \frac{p}{p+1}\sqrt{p} \cdot \frac{\ln(4)+2}{\ln(4p)+2} \qquad (p \geq 3,\,\nu_p(d) = 0).
    \end{equation}
    If \(\nu_p(d) =1\), then \(d \geq p\), so by \eqref{eq:PhiI.p>=3,nup(d)>0},
    \begin{equation} \label{eq:Phi2.nu=1}
        \frac{\tilde{\Phi}_2^\I(pd)}{\tilde{\Phi}_2^\I(d)} \geq \frac{p+1}{2p} \sqrt{p} \cdot \frac{\ln(4p)+2}{\ln(4p^2)+2} \qquad (p \geq 3,\,\nu_p(d) = 1).
    \end{equation}
    If \(\nu_p(d) \geq 2\), then \(d \geq p^2\) and \(\mu(2,p,d) = \mu(2,p,pd) = 0\), hence by \eqref{eq:PhiI.p>=3,nu>0(a)},
    \begin{equation} \label{eq:Phi2.nu>=2}
        \frac{\tilde{\Phi}_2^\I(pd)}{\tilde{\Phi}_2^\I(d)} \geq \sqrt{p} \cdot \frac{\ln(4p^2)+2}{\ln(4 p^3)+2} \qquad (p \geq3, \, \nu_p(d) \geq 2).
    \end{equation}
    Using \eqref{eq:1>=x/x+1>=x+1/2x} and
    \[
        \frac{\ln(4)+2}{\ln(4x)+2} \leq \frac{\ln(4x)+2}{\ln(4x^2)+2} \leq \frac{\ln(4x^2)+2}{\ln(4x^3)+2} \qquad (x \geq 1),
    \]
    we find 
    \begin{equation*} \label{eq:Phi2.p>=3}
        \frac{\tilde{\Phi}_2^\I(pd)}{\tilde{\Phi}_2^\I(d)} \geq \frac{p+1}{2p} \sqrt{p} \cdot \frac{\ln(4)+2}{\ln(4p)+2} \qquad (p \geq 3).
    \end{equation*}
    The function appearing on the right side is increasing, yielding
    \[
        \varphi(p) \geq \varphi(11) = \frac{12}{22} \sqrt{11}\cdot \frac{\ln(4)+2}{\ln(44)+2} \approx 1,059 > 1 \qquad (p \geq 11). 
    \]
    Define \(\varphi(p)\) as in the statement. It follows from \eqref{eq:Phi2.nu>=2} that
    \[
        \frac{\tilde{\Phi}_2^\I(3d)}{\tilde{\Phi}_2^\I(d)} \geq \varphi(3) \approx 1.4472 > 1 \qquad (\nu_3(d) \geq 2).
    \]
    It follows from \eqref{eq:Phi2.nu=0}, \eqref{eq:Phi2.nu=1}, \eqref{eq:Phi2.nu>=2} that
    \begin{align*}
        \frac{\tilde{\Phi}_2^\I(5d)}{\tilde{\Phi}_2^\I(d)} & \geq \min\bigg(\frac{5}{6} \sqrt{5} \cdot \frac{\ln(4)+2}{\ln(20)+2}, \varphi(5), \sqrt{5} \cdot \frac{\ln(100)+2}{\ln(500)+2} \bigg) = \varphi(5) \approx 1.0147 >1, \\
        \frac{\tilde{\Phi}_2^\I(7d)}{\tilde{\Phi}_2^\I(d)} & \geq \min\bigg(\frac{7}{8} \sqrt{7} \cdot \frac{\ln(4)+2}{\ln(28)+2}, \varphi(7), \sqrt{7} \cdot \frac{\ln(196)+2}{\ln(1372)+2} \bigg) = \varphi(7) \approx 1.1076 >1.
    \end{align*}
    By \Cref{lem:Phi_new}, \(\tilde{D}_2^\I\) is finite.

    Finally, we consider \(\Phi(d) = \tilde{\Phi}_2^\II(d)\). By definition of \(\tilde{\Phi}_2^\II\), we have  
    \[
        \frac{\tilde{\Phi}_2^\II(pd)}{\tilde{\Phi}_2^\II(d)} = \frac{\tilde{\Phi}_2^\I(pd)}{\tilde{\Phi}_2^\I(d)} \cdot \frac{2^{\max(0, 2 - \nu_2(pd))}}{2^{\max(0, 2 - \nu_2(d))}}.
    \]  
    If \(\nu_2(d) \geq 2\), then \(\nu_2(2d) \geq 3\), so we have  
    \[
        \frac{\tilde{\Phi}_2^\II(pd)}{\tilde{\Phi}_2^\II(d)} = \frac{\tilde{\Phi}_2^\I(pd)}{\tilde{\Phi}_2^\I(d)} \qquad (\nu_2(d) \geq 2),
    \]  
    and the same inequalities as above apply.
    By \Cref{lem:Phi_new}, \(\tilde{D}_2^\II\) is finite.
\end{proof}

\begin{lemma} \label{lem:estimate.zeta}
The function \(\zeta(2s)/\zeta(s)\) is strictly increasing for \(s > 1\). 
In particular, 
\[
    \frac{\zeta(2s)}{\zeta(s)} \geq \frac{\pi^2}{15} \qquad (s \geq 2).
\]
\end{lemma}

\begin{proof}
    Using Dirichlet convolution, it is easy to prove (see, e.g., \cite[formula (4) on p. 145]{hua}) that
    \[
        \frac{\zeta(s)}{\zeta(2s)} = \sum_{n = 1}^\infty \frac{|\mu(n)|}{n^s} \qquad (s > 1),
    \]
    where \(\mu\) is the Möbius function, defined as \(\mu(n) = (-1)^k\) if \(n\) is the product of \(k\) different primes, and \(\mu(n) = 0\) otherwise. Since \(n^s < n^t\) for \(s < t\), the result follows.
    
    For the last estimate, we use the well-known values \(\zeta(2) = \frac{\pi^2}{6}\) and \(\zeta(4) = \frac{\pi^4}{90}\). 
\end{proof}

\begin{lemma} \label{lem:estimate.Gamma}
    There exists \(x_0 \in \IR\) with \(\frac 72 < x_0 < 4\), such that the function~\(\pi^{-x}\Gamma(x)\) is increasing for \(x > x_0\).
\end{lemma}
\begin{proof}
    We denote by \(\psi(x) \coloneqq \Gamma'(x)/\Gamma(x)\) the digamma function. We refer to \cite[§12.3]{Whittaker.Watson:course.modern.analysis} for basic properties of \(\psi\). The integral representation
    \[
        \psi(x) = \int _{0}^{\infty }\left({\frac {e^{-t}}{t}}-{\frac {e^{-xt}}{1-e^{-t}}}\right)\,dt \qquad (x>0).
    \]
    implies that \(\psi\) is strictly increasing for \(x > 0\).
    Furthermore, we have
    \[
        (\pi^{-x}\Gamma(x))' = \pi^{-x}\Gamma(x) (\psi(x) - \ln(\pi)).
    \]
    Hence, the derivative of \(\pi^{-x}\Gamma(x)\) is positive for \(x > x_0\), where \(x_0 \in \IR\) is the unique positive value satisfying \(\psi(x_0) = \ln(\pi)\). 
    
    From the recurrence relation \(\psi(x+1) = \psi(x) + 1/x\), it follows that for any positive integer \(n\), we have
    \[
        \psi(n) = -\gamma + \sum_{k = 1}^{n-1} \frac 1 k, \quad \text{and} \quad \psi\big( n + \tfrac 12 \big) = -\gamma - 2\ln(2) + \sum_{k = 1}^n \frac{2}{2k-1},
    \]
    where \(\gamma = -\psi(1) \approx 0.577\) denotes the Euler--Mascheroni constant. We compute explicitly the following values:
    \[
        \psi\big(3 + \tfrac 12\big) \approx 1.103, \quad \ln(\pi) \approx 1.144, \quad \psi(4) \approx 1.256,
    \]
    which imply that \(3 + \frac 12 < x_0 < 4\).
\end{proof}

\begin{proposition} \label{prop:c_n>2}
    For every \(n \geq 19\), the inequality 
    \[
        c_n > 2
    \]
    holds, where \(c_n\) denotes the sequence defined in \Cref{thm:bound.rank.>2}.
\end{proposition}

\begin{proof}
Define the sequence 
\[
    b_n \coloneqq \pi^{-\frac{1}{2}(n+1)} \Gamma\big(\tfrac{1}{2}(n+1)\big) \frac{\zeta(2s)}{\zeta(s)} = \frac{c_n}{1 + 2^{-s}}.
\]
Clearly, \(c_n > b_n\) for all \(n\). By \Cref{lem:estimate.zeta,lem:estimate.Gamma}, the sequence \(b_n\) is monotonic increasing for \(n \geq 7\). Therefore, for \(n \geq 19\),
\[
    c_n > b_n \geq b_{19}.
\]

Next, compute \(b_{19}\):
\[
    b_{19} = \pi^{-10} \Gamma(10) \frac{\zeta(20)}{\zeta(10)}.
\]
Using \(\Gamma(10) = 9!\) and the estimate provided in \Cref{lem:estimate.zeta}, we find
\[
    b_{19} \geq \pi^{-10} \cdot 9! \cdot \frac{\zeta(4)}{\zeta(2)} \approx 2.5496 > 2.
\]
Thus, for all \(n \geq 19\), we have \(c_n > 2\), as required. \qedhere
\end{proof}

\begin{theorem} \label{thm:slender}
There exist only finitely many slender genera of rank \(n \geq 2\). 
\end{theorem}

In the following proof, we show that there are no slender genera of rank \(n \geq 19\).
In §\ref{subsec:enumeration}, we compute that there are no slender genera of rank \(n = 17,18\) either. 

\begin{proof}[Proof of \Cref{thm:slender}]
Recall that for any \(n\) and \(d\), there are only finitely many genera of rank \(n\) with determinant \(d\) (see, e.g., \cite[(20.2)]{kneser}). By \Cref{cor:DIn,cor:DI2}, the determinants of all slender genera of rank \(n \geq 2\) are contained within the sets \(D_n^\I\) and \(D_n^\II\). By \Cref{prop:long}, these sets are finite. Therefore, it suffices to show that \(D_n^\I\) and \(D_n^\II\) are empty for sufficiently large \(n\).

Define \(b = 1/c_n\). Consider first the set \(D_n^\I\). From \Cref{prop:long}, \(\Phi_n^\I(d)\) achieves its minimum value at \(d = 1\). Consequently, \(D_n^\I = \emptyset\) if 
\[
     \Phi_n^\I(1) = \frac{1}{2} > b,
\]
which holds for all \(n \geq 19\) by \Cref{prop:c_n>2}.

Similarly, for the set \(D_n^\II\), we have \(D_n^\II = \emptyset\) for all \(n \geq 19\), as
\[
     \Phi_n^\II(2^n) = \frac{\sqrt{2^n}}{8} > \frac{1}{2} > b.
\]
This completes the proof.
\end{proof}

\begin{proof}[Proof of \Cref{thm:idoneal}]
    The only idoneal genus of rank \(1\) is the genus of \([1]\). 
    
    We already observed in \Cref{rmk:idoneal} that there exist a one-to-one correspondence between idoneal genera of rank~\(2\) and idoneal numbers. It is known (cf. \Cref{rmk:GRH}) that there exist \(65\), \(66\) or \(67\) idoneal numbers.

    By \Cref{lem:slender.parent}, every idoneal genus of rank \(r \geq 3\) has a slender parent of rank \(n \geq 2\). Thus, we conclude by \Cref{thm:slender}.
\end{proof}

\subsection{Enumeration} \label{subsec:enumeration}

The proof of \Cref{thm:slender} suggests an algorithm to explicitly enumerate all slender genera of rank \(n \geq 3\), which we outline here.

\begin{lemma} \label{lem:possible_det}
    For any \(n \geq 2\), there exist algorithms to compute the sets \(D_n^{\I}\) and \(D_n^{\II}\), as defined in \Cref{cor:DIn,cor:DI2}.
\end{lemma}

\begin{proof}
    Using the algorithm described in \Cref{lem:Phi_new} along with the data provided in \Cref{prop:long}, we can compute the sets \(D_n^{\I}\) and \(D_n^{\II}\) for all \(n \geq 3\), as well as the auxiliary sets \(\tilde{D}_2^{\I}\) and \(\tilde{D}_2^{\II}\). 

    To determine \(D_2^{\I}\) and \(D_2^{\II}\), it suffices to examine all elements in the larger sets \(\tilde{D}_2^{\I}\) and \(\tilde{D}_2^{\II}\), selecting those that satisfy the extra condition. Hence, the sets \(D_n^{\I}\) and \(D_n^{\II}\) are computable for all \(n \geq 2\).
    
    In order to obtain proven results, the computations of \(c_n\), \(c_n\Phi^\I_n\) and \(c_n\Phi^\II_n\) are carried out using interval arithmetic with exact error bounds. Since we only need to prove the inequality $c_n\Phi^{\I / \II}_n\leq 1$, this is enough.  
    We used the library \verb+arb+  \cite{Johansson2017arb} via \texttt{sageMath} \cite{sage}.
    
    For \(2 \leq n \leq 18\), the output is summarized in \Cref{tab:possible.determinants}, where we also listed approximations of \(c_n\). 
    
    The computation of $D_2^\I$ and $D_2^\II$ is rather expensive and took about 50 days of CPU time. We employed the method \verb+quadratic_L_function__exact+ implemented by John Hanke \cite{hanke} in SageMath \cite{sage}, which uses the fact that 
    \(\zeta_D(1)\sqrt{d}\pi^{-1} \in \QQ\) (see, e.g., \cite[Chapter 12, Theorem 10.1]{hua}).
    The intermediate results are
    \begin{align*}
        |\tilde D_2^\I| & = 180791, & \max \tilde D_2^\I &= 45090045,  \\
        |\tilde D_2^\II| & = 1259915, & \max \tilde D_n^\II& =3607203600. \qedhere
    \end{align*}
\end{proof}

\begin{table}
\footnotesize
\caption{Sets of possible determinants of (odd/even) slender genera of rank \(n\), computed with the given approximation \(c\) of \(c_n\). }
\label{tab:possible.determinants}
\begin{tabular}{lrrrrrrrrr}
    \toprule
    \(n\)                 & 2 & 3      & 4      & 5      & 6      & 7      & 8      & 9     & 10    \\
    \(c\) &    0.1591      & 0.0833 & 0.0625 & 0.0614 & 0.0575 & 0.0607 & 0.0664 & 0.078 & 0.096 \\
    \midrule
    \(|D_n^\I|\) & 7939    & 960    & 1257   & 1210   & 1279   & 1126   & 923    & 665   & 435   \\
    \(\max D_n^\I\) & 475475 & 6750   & 5625   & 4131   & 4374   & 2187   & 1458  & 729   & 441   \\
    \(|D_n^\II|\) & 103750   & 2769   & 1846   & 900    & 476    & 209    & 85     & 31    & 9     \\
    \(\max D_n^\II\)&20263320  & 189000 & 90000  & 62208  & 46656  & 20736  & 14592  & 10752 & 6144  \\
    \midrule
    \(n\)                 & 11     & 12     & 13     & 14     & 15    & 16     & 17     & 18     \\
    \(c\)          & 0.1246 & 0.1687 & 0.2382 & 0.3493 & 0.531 & 0.8343 & 1.3525 & 2.2577 \\
    \midrule
    \(|D_n^\I|\)     & 258    & 140    & 70     & 32     & 14    & 5      & 2      & 0      \\
    \(\max D_n^\I\)  & 258    & 140    & 70     & 32     & 14    & 5      & 2      & --     \\
    \(|D_n^\II|\)    & 3      & 0      & 0      & 0      & 0     & 0      & 0      & 0      \\
    \(\max D_n^\II\) & 4096   & --     & --     & --     & --    & --     & --     & --     \\
    \bottomrule
\end{tabular}
\end{table} 

\begin{proposition} \label{prop:algorithm.slender}
    For every \(n \geq 2\), an algorithm exists to generate a complete list of all slender genera of rank~\(n\).
\end{proposition}

\begin{proof}
    Let \(S\) be the output of the following algorithm:
    \begin{myalgorithm}\caption{Slender genera}\label{alg:slender}
    \begin{algorithmic}[1]
         \REQUIRE an integer \(n \geq 2\), the sets \(D_n^{\I}\) and \(D_n^\II\).
         \ENSURE the list of all slender genera \(g\) of rank~\(n\).
         \STATE \(S \longleftarrow \emptyset\)
        \STATE \(G^{\I} \longleftarrow\) list of all odd genera of rank~\(n\) and determinant \(d \in D_n^{\I}\) \label{step:list_detI}
         \STATE \(G^{\II} \longleftarrow\) list of all even genera of rank~\(n\) and determinant \(d \in D_n^{\II}\) \label{step:list_detII}
        \FORALL{\(g \in G^{\I} \cup G^{\II}\)}
            \STATE \(g' \longleftarrow\) the child of \(g\) (\Cref{def:parent}) \label{step:child}
            \IF{\(m(g') \leq m(g)\)} \label{line:def.slender}
                \STATE \(S \longleftarrow S \cup \{g\}\)
            \ENDIF
        \ENDFOR
         \RETURN \(S\)
    \end{algorithmic}
    \end{myalgorithm}

    By \Cref{cor:DIn,cor:DI2}, the determinant \(d\) of a slender genus of rank \(n\), odd or even respectively, belongs to \(D_n^\I\) or \(D_n^\II\). 

    \sloppy
    Lines \ref{step:list_detI} and \ref{step:list_detII} can be executed with the \verb+genera+ function provided in \verb+sage.quadratic_forms.genera.genus,+ implemented by Brandhorst \cite{brandhorst_discriminant_forms}.
    Thus, if \(g\) is a slender genus, then \(g \in G^\I \cup G^\II\). 

    A genus is represented in terms of its Conway-Sloane genus symbol \cite[Chapter 15 \S 7.8]{splag}.
    Using the function \verb+direct_sum+ of the class \verb+GenusSymbol_global_ring+, it is not necessary to compute a representative of \(g\) (which would take up much time) to compute \(g'\) in line~\ref{step:child}. Likewise, the mass can be computed directly from the genus symbol. The computations for \(n = 2\) took about 25 CPU days. For \(n > 2\) computations terminated in a matter of hours.

    In line \ref{line:def.slender}, the definition of slender genus is applied (\Cref{def:slender}). Thus, \(S\) is a list of all slender genera of rank \(n\). 
    
    The output is summarized in \Cref{tab:summary.slender}.     
\end{proof}

\begin{table}
\footnotesize
\caption{Sets \(S^\I_n\), \(S^\II_n\) of odd (resp. even) slender genera of rank \(n\).}
\label{tab:summary.slender}
\begin{tabular}{lrrrrrrrrr}
    \toprule
    \(n\)         & 2     & 3    & 4    & 5    & 6    & 7    & 8    & 9    & 10  \\
    \midrule 
    \(|S_n^\I|\)  &  306   & 782  & 1406 & 2069 & 2374 & 2288 & 1804 & 1230 & 712 \\
    \(\max \det S_n^\I\)   & 2700 & 2268 & 2592 & 1728 & 1152 & 1024 & 768  & 320  & 128 \\
    \(|S_n^\II|\) &  244   & 195  & 130  & 65   & 22   & 8    & 4    & 0    & 0     \\
    \(\max \det S_n^\II\) & 10800& 6912 & 4096 & 4096 & 1024 & 1024 & 1024 & --   & -- \\
    \midrule 
    \(n\)       & 11  & 12  & 13 & 14 & 15 & 16 & 17 \\
    \midrule 
    \(|S_n^\I|\)  & 361 & 157 & 61 & 19 & 4  & 1  & 0  \\
    \(\max \det S_n^\I\)  & 68  & 36  & 17 & 8  & 3  & 1  & -- \\
    \(|S_n^\II|\) & 0   & 0   & 0  & 0  & 0  & 0  & 0  \\
    \(\max \det S_n^\II\) & --  & --  & -- & -- & -- & -- & -- \\
    \bottomrule
\end{tabular}
\end{table} 

\begin{theorem} \label{thm:algorithm.idoneal}
    For all \(n \geq 3\), there exists an algorithm to compute the list of all idoneal genera of rank \(n\).
\end{theorem}
\begin{proof}
    Let \(I\) be the output of the following algorithm:
    \begin{myalgorithm}\caption{Idoneal genera}\label{alg:idoneal_genera}
    \begin{algorithmic}[1]
        \REQUIRE an integer \(n \geq 3\), the list \(S\) of all slender genera of rank \(n-1\).
        \ENSURE the list of all idoneal genera of rank~\(n\).
        \STATE \(Q,T,I \longleftarrow \emptyset\)
        \FORALL{\(f_1 \in S\)} \label{line:genera_with_slender_parents}
            \STATE \(g \longleftarrow\) child of \(f_1\)
            \IF{\(g\) has another parent \(f_2\)} \label{line:if_g_has_another_parent}
                \STATE \(Q \longleftarrow Q \cup \{(g, \{f_1,f_2\})\}\)
            \ELSE
                \STATE \(Q \longleftarrow Q \cup \{(g,\{f_1\})\}\) 
            \ENDIF
        \ENDFOR
        \FORALL{\((g,F) \in Q\)} \label{line:sieve}
            \IF{\(2m(g) \leq \sum_{f \in F} m(f)\)} \label{eq:sieve}
                \STATE \(T \longleftarrow T \cup \{(g,F)\}\)
            \ENDIF
        \ENDFOR
        \FORALL{\((g,F) \in T\)} \label{step:loop_representatives}
            \STATE \(R \longleftarrow\) list of representatives of all lattices in the genera contained in \(F\) \label{step:representatives}
            \STATE \(m \longleftarrow \displaystyle \sum_{L \in R} \frac{1}{|{\Aut(L \oplus [1])}|}\) \label{step:m}
            \IF{\(m = m(g)\)}
                \STATE \(I \longleftarrow I \cup \{g\}\)
            \ENDIF
        \ENDFOR
        \RETURN \(I\)
    \end{algorithmic}
    \end{myalgorithm}

    We now explain the steps of the algorithm and justify its correctness.

    By \Cref{lem:slender.parent}, every idoneal genus of rank \(n\) has at least one slender parent of rank \(n-1\). 

    In the loop starting at line \ref{line:genera_with_slender_parents}, we build the set \(Q\), consisting of all pairs \((g, F)\), where \(g\) is a genus of rank \(n\) with at least one slender parent, and \(F\) is the set of its parents. 
    By \Cref{lem:at.most.two.parents}, the set \(F\) contains one or two genera. 
    The condition on line \ref{line:if_g_has_another_parent} is straightforward to check, as \(f_2\), if it exists, must have the same rank, determinant, and bilinear discriminant form as \(f_1\), but opposite parity (see \Cref{lem:at.most.two.parents}).

    The loop starting at line \ref{line:sieve} filters the pairs \((g, F) \in Q\) by applying the necessary condition on the masses, as given in \Cref{lem:slender.parent}. This step is fast since it relies on mass calculations for which we used the code written in \verb+sageMath+ by John Hanke~\cite{hanke}.

    Finally, in the loop starting at line \ref{step:loop_representatives}, we compute the list \(R\) of representatives for all lattices in the genera in \(F\). If \(g\) is idoneal, all its representatives are of the form \(L \oplus [1]\), where \(L\) belongs to some genus of \(F\). Thus, the value \(m\) computed in line \ref{step:m} will equal \(m(g)\). Otherwise, \(m < m(g)\), and \(g\) is excluded from the output. Therefore, \(I\) correctly contains all idoneal genera of rank \(n\).

    The process is summarized in \Cref{tab:summary.idoneal}, including intermediate steps. 
    
    Note that line \ref{step:representatives} is the most computationally intensive, as it involves enumerating representatives of a genus. 
    The computations took us roughly one month. We employed Kneser's neighbor method (\cite[§28]{kneser}, \cite{neighbor}, \cite{kirschmer}) for this step, using algorithms from \cite{fincke.pohst} and \cite{isometries} implemented in \verb+PARI+~\cite{pari} and, in difficult cases, \verb+Magma+~\cite{magma}. The most critical subtasks involve short vector enumeration, isometry testing, and computing orthogonal groups.
\end{proof}

\begin{remark}
    Of the 110 idoneal genera in rank \(3\), 15 are listed in Theorem 48 of \cite{Kani:idoneal.numbers}. Our findings align precisely with the results of that theorem in the specific case considered there.
\end{remark}

\begin{table}
\footnotesize
\caption{Sets \(Q_n\) of genera of rank \(n\) with at least one slender parent, subsets \(T_n \subseteq Q_n\) of genera satisfying the equation in Line~\ref{eq:sieve} of \Cref{alg:idoneal_genera} , and subsets \(I_n \subseteq T_n\) of idoneal genera of rank \(n\). }
\label{tab:summary.idoneal}
\begin{tabular}{lrrrrrrrrrrrr}
    \toprule
    \(n\)             & 3    & 4    & 5    & 6    & 7    & 8    & 9    & 10  \\
    \midrule 
    \(|Q_n|\)         & 543   & 965  & 1522 & 2125 & 2394 & 2295 & 1807 & 1230 \\
    \(\max \det Q_n\) & 10800   & 6912 & 4096 & 4096 & 1152 & 1024 & 1024 & 320  \\
    \(|T_n|\)         & 142   & 196  & 297  & 393  & 425  & 400  & 300  & 199  \\
    \(\max \det T_n\) & 1200   & 768  & 1024 & 256  & 256  & 256  & 256  & 64   \\
    \(|I_n|\)         & 110& 122  & 107 & 76    & 47 & 24 & 13 & 6  \\
    \(\max \det I_n\) & 1200& 768 & 1024 & 256   & 256& 256& 256& 8   \\
    \midrule 
    \(n\)             & 11  & 12  & 13  & 14 & 15 & 16 & 17 \\
    \midrule 
    \(|Q_n|\)         & 712 & 361 & 157 & 61 & 19 & 4  & 1  \\
    \(\max \det Q_n\) & 128 & 68  & 36  & 17 & 8  & 3  & 1  \\
    \(|T_n|\)         & 106 & 55  & 20  & 8  & 2  & 0  & 0  \\
    \(\max \det T_n\) & 32  & 16  & 8   & 4  & 2  & -- & -- \\
    \(|I_n|\)         & 4 & 1 & 1 & 0  & 0  & 0  & 0  \\
    \(\max \det I_n\) & 4  & 2 & 1  & -- & -- & -- & -- \\
    \bottomrule
\end{tabular}
\end{table}

\section{K3 surfaces covering an Enriques surface} \label{sec:K3-covering-Enriques}

The aim of this section is to prove \Cref{thm:transcendental}, which is done in §\ref{subsec:proof-transc-thm}. 
First, we review some fundamental properties of finite quadratic forms and Nikulin's theory of discriminant forms in §§\ref{subsec:finite-quadratic-forms}--\ref{subsec:nikulin}. 
Next, the possible shapes of the discriminant form of lattices with signature \((2, \lambda - 2)\) embedding in \(\Lminus\) are determined in \Cref{sec:T-embedding-in-Lminus} and listed in \Cref{tab:T-embedding-in-Lminus}. 
Finally, we identify all co-idoneal lattices in §\ref{subsec:co-idoneal}.

\subsection{Finite quadratic forms} \label{subsec:finite-quadratic-forms}
Let \(A\) be a finite commutative group. The bilinear form associated to a finite quadratic form \(q\colon A \rightarrow \IQ/2\IZ\) is given by
\[
q^\flat\colon A \times A \rightarrow \IQ/\IZ, \quad q^\flat(\alpha,\beta)\coloneqq \tfrac{1}{2} (q(\alpha+\beta)-q(\alpha)-q(\beta)).
\] 
A torsion quadratic form \(q\) is \emph{nondegenerate} if the homomorphism of groups \(A \rightarrow \Hom(A,\IQ/\IZ)\) given by \(\alpha \mapsto (\beta \mapsto q^\flat(\alpha, \beta))\) is an isomorphism.
If \(H \subset A\) is a subgroup, then \(q|H\) denotes the restriction of \(q\) to \(H\).

For \(k \in \IN\), we denote by \(\bu_k\), \(\bv_k\) the finite quadratic forms with the underlying group 
\(\IZ/2^k\IZ e_1 \oplus \IZ/2^k\IZ e_2\) and quadratic form 
\[
    \bu_k(x e_1 + y e_2) = \frac{1}{2^{k-1}} xy, \qquad \bv_k(x e_1 + y e_2) = \frac{1}{2^{k-1}} (x^2 + xy + y^2).
\]
For \(k \in \IN\) and \(\varepsilon \in \{1,3,5,7\}\), we denote by \(\bw^\varepsilon_{2,k}\) the finite quadratic forms with underlying group \(\IZ/2^k\IZ e\) and
\[
\bw^\varepsilon_{2,k}(xe) = \frac{\varepsilon}{2^{k}} x^2.
\]
Note that \(\bw^1_{2,1} \cong \bw^5_{2,1}\) and \(\bw^3_{2,1} \cong \bw^7_{2,1}\). For an odd prime~\(p\), \(k \in \IN\) 
and \(\varepsilon \in \{\pm 1\}\), we denote by \(\bw^\varepsilon_{p,k}\) the finite quadratic forms with underlying group \(\IZ/p^k\IZ e\) and 
\[
    \bw^\varepsilon_{p,k}(xe) = \frac{a}{p^{k}} x^2,
\] 
where \(a\) is an integer with Legendre symbol \(\Big(\displaystyle \frac{2a}{p}\Big) = \varepsilon\). 

Moreover, the degenerate quadratic form \(q\colon A \rightarrow \IZ/2\IZ\) with \(A \cong \IZ/2\IZ\) taking the values \(0 \mod 2 \IZ\) (resp. \(1 \mod 2\IZ\)) on the nontrivial element is denoted by~\(\langle 0 \rangle\) (resp. \(\langle 1 \rangle\)).
If \(A\) is \(2\)-elementary, then the short exact sequence 
\[
0 \rightarrow q^\perp \rightarrow q \rightarrow q/q^\perp \rightarrow 0
\]
splits.
Hence, \(q\) can be written as the direct sum of copies of \(\bu_1, \bv_1, \bw^1_{2,1}, \bw^3_{2,1}, \langle 0 \rangle, \langle 1 \rangle\).
In this paper, a (possibly degenerate) finite quadratic form \(q\) is called \emph{odd} if \(q \cong \bw^\varepsilon_{2,1} \oplus q'\) for some~\(\varepsilon\) and some finite quadratic form~\(q'\), otherwise it is called \emph{even}.
From now on, unless explicitly stated, we assume a quadratic form to be nondegenerate.

The abelian group \(A = L^\vee/L\), where 
\[
    L^\vee \coloneqq \{v \in L \otimes \IQ \mid \text{\(b(v,w) \in \IZ\) for all \(w \in L\)}\},
\]
has order \(\det(L)\).
If \(L\) is an even lattice, the finite quadratic form \(q(L)\colon A \rightarrow \IQ/2\IZ\) induced by the linear extension of \(b\) to~\(\IQ\) is called the \emph{discriminant (quadratic) form} of~\(L\). 
The \emph{signature} of a finite quadratic form \(q\) is defined as \(\sign q = s_+ - s_- \mod 8\), where \((s_+,s_-)\) is the signature of any even lattice \(L\) such that \(q(L) \cong q\). 
Explicit formulas for the signature of elementary quadratic forms are provided in \cite[Proposition 1.11.2]{Nikulin:int.sym.bilinear.forms}). 
Given a prime number \(p\), \(q_p\) denotes the restriction of \(q\) to the \(p\)-Sylow subgroup. 
We put \(\ell_p(q) = \ell(q_p)\), where \(\ell\) denotes the length of a finite abelian group.

\begin{lemma} \label{lem:parity-l2}
For any finite quadratic form \(q\), \(\ell_2(q) \equiv \sign q \mod 2\).
\end{lemma}
\proof 
This can be checked for all \(\bu_k\), \(\bv_k\), \(\bw^\varepsilon_{p,k}\). The claim follows by linearity, as each finite quadratic form is isomorphic to the direct sum of such forms (see, e.g., \cite[Proposition~1.8.1]{Nikulin:int.sym.bilinear.forms}). 
\endproof

\begin{lemma} \label{lem:2-elementary-even-subgroup}
Any finite quadratic form \(q\) contains a \(2\)-elementary subgroup~\(H\) of length \(\ell_2(q) - 1\) such that \(q|H\) is even. A finite quadratic form \(q\) contains a \(2\)-elementary subgroup \(H\) of length \(\ell_2(q)\) such that \(q|H\) is even if and only if \(q\) is even.
\end{lemma}
\proof
Let \(G\) be the underlying group of \(q\).
The subgroup \(H = \{\alpha \in G: 2\alpha = 0\}\) is \(2\)-elementary of length \(\ell_2(q)\). The (possibly degenerate) quadratic form \(q|H\) is even if and only if \(q\) is. This proves the second statement.

Write \(q|H\) as the direct sum of copies of \(\bu_1, \bv_1, \bw^1_{2,1}, \bw^3_{2,1}, \langle 0 \rangle, \langle 1 \rangle\).
By \cite[Proposition~1.8.2]{Nikulin:int.sym.bilinear.forms}, we can suppose that there are at most two copies of \(\bw^\varepsilon_{2,1}\). 
If there is only one copy of \(\bw^\varepsilon_{2,1}\), then the underlying subgroup \(H'\) to the rest of the direct sum has length \(\ell_2(q)-1\) and \(q|H'\) is even. 
If there are two copies of \(\bw^\varepsilon_{2,1}\), and \(H''\) is the underlying subgroup of length \(\ell_2(q) -2\) to the rest of the direct sum, we observe that \(\bw^1_{2,1} \oplus \bw^1_{2,1}\) and \(\bw^3_{2,1} \oplus \bw^3_{2,1}\) contain a copy of \(\langle 1 \rangle\), while \(\bw^1_{2,1} \oplus \bw^3_{2,1}\) contains a copy of \(\langle 0 \rangle\), so we take either \(\langle 1 \rangle \oplus H''\) or \(\langle 0 \rangle \oplus H''\), and we conclude. 
\endproof 

The following lemma uses similar ideas as the previous ones. We leave the details to the reader.

\begin{lemma} \label{lem:subgroups-u1^n}
Let \(q = n\bu_1\). If \(H\) is a subgroup of \(q\) and \(m = \max\{ 0, \ell_2(H) - n\}\), then there exists a (possibly degenerate) finite quadratic form \(q'\) such that 
\[
\pushQED{\qed} 
     q|H \cong m\bu_1 \oplus q'. \qedhere
\popQED
\]
\end{lemma}

\subsection{Nikulin's theory of discriminant forms} \label{subsec:nikulin}

We recall here some basic results on the theory of discriminant forms as developed by Nikulin~\cite{Nikulin:int.sym.bilinear.forms}.

\begin{theorem}[Nikulin {\cite[Theorem 1.9.1]{Nikulin:int.sym.bilinear.forms}}] 
For each finite quadratic form~\(q\) and prime number \(p\) there exists a unique \(p\)-adic lattice \(K_p(q)\) of rank \(\ell_p(q)\) whose discriminant form is isomorphic to \(q_p\), except in the case when \(p = 2\) and \(q\) is odd. \qed
\end{theorem}

\newcommand{\condA}[1]{\hyperlink{cond:A}{\mathrm A}(#1)}
\newcommand{\condB}[1]{\hyperlink{cond:B}{\mathrm B}(#1)}
\newcommand{\condC}[1]{\hyperlink{cond:C}{\mathrm C}(#1)}

We introduce the following conditions (depending on \(s,s' \in \IZ\)) on a finite quadratic form \(q\). 

\vspace{3mm}
\begin{tabular}{ll}
    \hypertarget{cond:A}{\({\mathrm A}(s)\)}: &  \(\sign q \equiv s \mod 8\). \\[3mm]
    \hypertarget{cond:B}{\({\mathrm B}(s, s')\)}: & \parbox{.79\textwidth}{for all primes \(p \neq 2\), \(\ell_p(q) \leq s + s'\); moreover, \\ \(|q| \equiv (-1)^{s'} \discr K_p(q) \mod (\IZ_p^\times)^2\) if \(\ell_p(q) = s + s'\).} \\[4mm]
    \hypertarget{cond:C}{\({\mathrm C}(s)\)}: & \parbox{.79\textwidth}{\(\ell_2(q) \leq s\); moreover, \(|q| \equiv \pm \discr K_2(q) \mod (\IZ_2^\times)^2\) \\ if \(\ell_2(q) = s\) and \(q\) is even.}
\end{tabular}
\vspace{3mm}

\begin{theorem}[Nikulin {\cite[Theorem 1.10.1]{Nikulin:int.sym.bilinear.forms}}] \label{thm:Nikulin-existence}
An even lattice of signature \((t_{+},t_{-})\), \(t_+,t_- \in \IZ_{\geq 0}\), and discriminant quadratic form \(q\) exists 
if and only if \(q\) satisfies conditions \(\condA{t_+ - t_-}\), \(\condB{t_+, t_-}\) and \(\condC{t_+ + t_-}\). \qed
\end{theorem}

\begin{theorem}[Nikulin {\cite[Theorem 1.14.2]{Nikulin:int.sym.bilinear.forms}}] \label{thm:Nikulin-uniqueness}
If \(T\) is an even, indefinite lattice satisfying the following conditions:
\begin{enumerate}[(a)]
 \item \(\rank T \ge \ell_p(q(T)) + 2\) for all \(p \neq 2\),
 \item if \(\rank T = \ell_2(q(T))\), then \(q(T) \cong \bu_1 \oplus q'\) or \(q(T) \cong \bv_1 \oplus q'\),
\end{enumerate}
then the genus of \(T\) contains only one class. \qed
\end{theorem}

Given a pair of nonnegative integers \((m_+, m_-)\) and a finite quadratic form~\(q\), Nikulin establishes a useful way to enumerate the set of primitive embeddings of a fixed lattice \(T\) into any even lattice belonging to the genus~\(g\) of signature \((m_+,m_-)\) and discriminant form~\(q\). 
Nikulin's proposition is hindered though by a poor translation. In the original Russian text \cite{Nikulin:int.sym.bilinear.forms_Russian}, the term ``четные решетки'' has been incorrectly rendered in English as the singular phrase ``an even lattice.'' However, the correct translation is ``even lattices'' (plural). 
For the sake of completeness, we restate the proposition correctly:

\begin{proposition}[{\cite[Proposition 1.15.1]{Nikulin:int.sym.bilinear.forms}}] \label{prop:Nikulin1.15.1}
Let \(T\) be an even lattice of signature \((t_+,t_-)\) and \(g\) an even genus of signature \((m_+,m_-)\) and discriminant quadratic form \(q\). 
Then, for each lattice~\(S\), there exists a primitive embedding \(T \hookrightarrow \Lambda\) with \(T^\perp \cong S\) and \(\Lambda\) a lattice of genus \(g\) if and only if \(\sign S = (m_+ - t_+, m_- - t_-)\) and
there exist subgroups \(H \subset q\) and \(K \subset q(T)\), and an isomorphism of quadratic forms \(\gamma\colon q|H \rightarrow q(T)|K\), whose graph is denoted by~\(\Gamma\), such that 
\[
\pushQED{\qed} 
 q(S) \cong (q \oplus (-q(T))) | \Gamma^\perp/\Gamma. \qedhere
\popQED
\]
\end{proposition}

The following proposition is a simplified version of \Cref{prop:Nikulin1.15.1} in the case that the genus~\(g\) contains only one class. In this paper, we will only use this version.

\begin{proposition} 
\label{prop:Nikulin1.15.1-discriminant-form-orthogonal}
Let \(T\) be an even lattice of signature \((t_+,t_-)\) and \(\Lambda\) be an even lattice of signature \((m_+,m_-)\) which is unique in its genus. 
Then, for each lattice~\(S\), there exists a primitive embedding \(T \hookrightarrow \Lambda\) with \(T^\perp \cong S\) if and only if \(\sign S = (m_+ - t_+, m_- - t_-)\) and
there exist subgroups \(H \subset q(\Lambda)\) and \(K \subset q(T)\), and an isomorphism of quadratic forms \(\gamma\colon q(\Lambda)|H \rightarrow q(T)|K\), whose graph is denoted by~\(\Gamma\subset q(\Lambda) \oplus -q(T)\), such that 
\[
\pushQED{\qed} 
 q(S) \cong (q(\Lambda) \oplus (-q(T))) | \Gamma^\perp/\Gamma. \qedhere
\popQED
\]
\end{proposition}

\subsection{Transcendental lattices embedding in \texorpdfstring{\(\Lambda^-\)}{Lambda minus}} \label{sec:T-embedding-in-Lminus}

Recall that we defined \(\Lminus \coloneqq \bU \oplus \bU(2) \oplus \bE_8(-2)\). It follows that \(\sign(\Lminus) = (2,10)\) and \(q(\Lminus) \cong 5\bu_1\).

In this section the possible discriminant forms of even lattices of signature \((2, \lambda - 2)\) which embed primitively into \(\Lminus\) are determined. Necessarily, \(2 \leq \lambda \leq 12\).

\begin{lemma} \label{lem:technical}
Let \(f \colon F \rightarrow \IQ/2\IZ\) and \(g \colon G \rightarrow \IQ/2\IZ\) be finite quadratic forms with \(f \cong n \bu_1\). 
Let \(H \subset F\) and \(K \subset G\) be subgroups and \(\gamma \colon f|H \rightarrow g|K\) an isometry. 
Let \(\Gamma\) be the graph of \(\gamma\) in \(F \oplus G\). Then
\[
 \ell_2(H)\bu_1 \oplus \left(f \oplus (-g) | \Gamma^\perp/\Gamma\right) \cong f \oplus (-g).
\]
In particular,
\[
 \ell_2(f) + \ell_2(g) = \ell_2(\Gamma^\perp/\Gamma) + 2\ell_2(H).
\]
\end{lemma}
\proof
Recall that \(f^\flat\) and \(g^\flat\) denote the bilinear forms induced by the quadratic forms \(f\) and \(g\), respectively.
Let \(C=\ker (f^\flat|H)\). 
The exact sequence 
\[
 0 \rightarrow C^\perp \rightarrow F \rightarrow \Hom(C, \IQ/\IZ) \rightarrow 0
\]
(induced by \(f^\flat\)) splits, because \(F\) is \(2\)-elementary by assumption.

Let \(s\colon \Hom(K,\IQ/\IZ) \rightarrow F\) be a section and \(C^\vee_s\) its image. 
Then, since \(C\) is a totally isotropic subspace with respect to \(f^\flat\), we infer that \(f^\flat|(C \oplus C_s^\vee) \cong \ell u_1^\flat\), with \(\ell = \ell_2(C)\).
By modifying the section \(s\), we may assume that \(f^\flat|C_s^\vee = 0\).
Consider the subgroups \(F_s' = (H \oplus C_s^\vee)^\perp \subset F\), \(D_s = (C \oplus C^\vee_s)^{\perp} \subset H\oplus C^\vee_s\) and \(G'_s = \gamma(D_s)^\perp \subset G\).
Putting \(f' = f|F'_s\), \(g' = g|G'_s\) and \(d = f|D_s\), we obtain
\begin{align*} 
 f^\flat &\cong f^\flat|D_s \oplus (C\oplus C_s^\vee) \oplus F'_s \cong d^\flat \oplus \ell u_1^\flat \oplus (f')^\flat, \\
 g^\flat & \cong g^\flat|\gamma(D_s) \oplus G'_s \cong d^\flat \oplus (g')^\flat.
\end{align*}
Let \(\varphi\colon G'_s \rightarrow C_s^\vee\) be defined by
\(f^\flat(\varphi(\beta), \alpha) = g^\flat(\beta, \gamma(\alpha))\)
for all \(\alpha \in C\subset H\) and \(\beta \in G'_s\).
Define \(\psi\colon G'_s \rightarrow \Gamma^\perp\) by \(\psi(\beta) = \varphi(\beta) + \beta\).
Since \(f^\flat|C^\vee_s = 0\), we have
\begin{equation} \label{eq:technical}
 f^\flat \oplus (-g^\flat)(\psi(\beta),\psi(\beta')) = f^\flat(\varphi(\beta),\varphi(\beta')) - g^\flat(\beta,\beta') = -g^\flat(\beta,\beta'). 
\end{equation}
It follows that \(f^\flat \oplus (-g^\flat)\) restricted to \(F'_s \oplus \psi(G'_s) \subset \Gamma^\perp\) is nondegenerate. Since the orders of \(\Gamma^\perp/\Gamma\) and \(F_s' \oplus \psi(G'_s)\) coincide, this shows that
\[
 (f \oplus (-g))^\flat|\Gamma^\perp/\Gamma \cong (f' \oplus (-g'))^\flat.
\]

We claim that we can choose the section \(s\) in such a way that the quadratic forms coincide.
Indeed, if \(f|C^\vee_s = 0\), then \eqref{eq:technical} also holds at the level of quadratic forms, so we can replace \(f^\flat,g^\flat\) by \(f,g\), respectively, and we are done (note that \(2d \cong (\ell_2(H)-\ell)\bu_1\)).

If \(C \oplus C_s^\vee = F\), the section \(s\) can be modified so that \(f|C_s^\vee = 0\), because \(f \cong n\bu_1\).
Otherwise, as \(f|(C\oplus C_s^\vee)^\perp\) is then even, nondegenerate and nonzero, 
there exists \(\alpha \in (C \oplus C_s^\vee)^\perp\subset C^\perp\) with \(f(\alpha) = 1\). 
Let \(\delta_1,\ldots,\delta_\ell\) be a basis of \(\Hom(C,\IQ/\IZ)\), 
so that \(C_s^\vee = \langle s(\delta_1),\ldots,s(\delta_\ell)\rangle\).
Define \(s'\) by \(s'(\delta_i) = s(\delta_i)\) if \(x(s(\delta_i)) = 0\) and \(s'(\delta_i) = s(\delta_i) + \alpha\) else.
By replacing \(s\) with \(s'\) we are again in the situation where \(x|C_s^\vee = 0\) and we conclude.
\endproof

\begin{proposition} \label{prop:T-embedding-in-Lminus}
An even lattice \(T\) of signature \((2, \lambda - 2)\) embeds primitively into~\(\Lminus\) if and only if \(q(T)\) is of the form given in \Cref{tab:T-embedding-in-Lminus} for some nondegenerate finite quadratic form \(q\) satisfying the given conditions. 
In that case, \(q(T^\perp)\) is isomorphic to the form given in the corresponding column of the table.
\end{proposition}

\begin{proof} 
If \(\lambda = 12\), then \(\rank T = \rank \Lminus\), so the claim is trivial. For the rest of the proof we will suppose that \(2 \leq \lambda \leq 11\).

Assume first that \(T\) embeds in \(\Lminus\) and let \(H \subset q(T)\) be the subset given by \Cref{prop:Nikulin1.15.1-discriminant-form-orthogonal}.
By the last equation in \Cref{lem:technical} we see that
\[
 \ell_2(q(T)) + 10 = \ell_2(q(T^\perp)) + 2\ell_2(H).
\]
Using \(\ell_2(q(T^\perp)) \leq \rank(T^\perp) = 12 - \lambda\) and \(\ell_2(H) \leq \ell_2(q(T))\), we see that
\[
 \lambda -2 = 10 - \rank(T^\perp)\leq 2\ell_2(H)-\ell_2(q(T)) \leq \ell_2(q(T)) \leq \rank(T) = \lambda.
\]
By \Cref{lem:parity-l2}, \(\ell_2(q(T))\) can only assume two values, namely \(\lambda - 2\) or \(\lambda\) . 

We also infer that if \(\ell_2(q(T)) = \lambda -2\), then \(\ell_2(H) = \lambda - 2\), 
whereas if \(\ell_2(q(T)) = \lambda\), then either \(\ell_2(H) = \lambda -1\) or \(\ell_2(H) = \lambda\) 
(if \(\lambda = 11\) only \(\ell_2(H) = \lambda -1\) is possible, as \(\ell_2(H) \leq \ell_2(q(\Lminus)) = 10\)).
Moreover, by \Cref{lem:2-elementary-even-subgroup}, \(q\) must be even whenever \(\ell_2(H) = \ell_2(q(T))\).

Summarizing, for each \(\lambda \in \{2, \ldots, 11 \}\) we have three cases, described in \Cref{tab:cases}, except for \(\lambda = 11\) where the last case does not occur.

\begin{table}
\caption{The three cases appearing in the proof of \Cref{prop:T-embedding-in-Lminus}.}
\label{tab:cases}
\begin{tabular}{lllll}
 \toprule
 \(\lambda_{\rm case}\) & \(\ell_2(q(T))\) & \(\ell_2(H)\) & \(\ell_2(q(T^\perp))\) & parity of \(q(T)\) \\
 \midrule
 \(\lambda_{\rm a}\) & \(\lambda - 2\) & \(\lambda - 2\) & \(12 - \lambda\) & even \\
 \(\lambda_{\rm b}\) & \(\lambda\) & \(\lambda - 1\) & \(12 - \lambda\) & even or odd \\
 \(\lambda_{\rm c}\) & \(\lambda\) & \(\lambda\) & \(10 - \lambda\) & even \\
 \bottomrule
\end{tabular}
\end{table}

Let \(m = \max\{ 0, \ell_2(H) - 5\}\). \Cref{lem:subgroups-u1^n} ensures the existence of a subgroup of \(H\) (hence of \(q(T)\)) isometric to \(m\bu_1\); therefore, \(q(T) \cong m\bu_1 \oplus q\). 
The form of \(q(T^\perp)\) is then given by \Cref{lem:technical} and we can apply \Cref{thm:Nikulin-existence} to find all necessary conditions on~\(q\). 
We only write those that are also sufficient: for instance, condition \(\condA{12 - \lambda}\) for \(q(T^\perp)\) is always equivalent to \(\condA{\lambda}\) for \(q(T)\); if \(2 \leq \lambda \leq 6\), condition \(\condB{2,\lambda -2}\) for \(q(T)\) automatically implies \(\condB{0,12-\lambda}\) for \(q(T^\perp)\). 

Conversely, if \(T\) is given as in one of the rows of \Cref{tab:T-embedding-in-Lminus}, then the conditions on~\(q\), together with \Cref{prop:Nikulin1.15.1-discriminant-form-orthogonal} and \Cref{lem:2-elementary-even-subgroup}, ensure the existence of a primitive embedding \(T \hookrightarrow \Lminus\) with \(q(T^\perp)\) of the given form.
\end{proof}

\begin{center}

\begin{table}
    \caption{Discriminant forms of lattices \(T\) of signature \((2,\lambda -2)\) embedding primitively into {\(\mathbf{\Lambda^-}\)} (see \Cref{prop:T-embedding-in-Lminus}).}
    \label{tab:T-embedding-in-Lminus}

    \begin{tabular}{lllll}
    \toprule
    \(\lambda_{\rm case}\) & \(q(T)\) & \(q(T^\perp)\) & \(\ell_2(q)\) & conditions on \(q\) \\
    \midrule

    \(2_{\rm a}\) & \(q\) & \(5\bu_1 \oplus q(-1)\)     & \(0\) & \(\condC{0}\) \\
    \(2_{\rm b}\) & \(q\) & \(4\bu_1 \oplus q(-1)\)    & \(2\) & -- \\
    \(2_{\rm c}\) & \(q\) & \(3\bu_1 \oplus q(-1)\)   & \(2\) & even \\

    \midrule
    
    \(3_{\rm a}\) & \(q\) & \(4\bu_1 \oplus q(-1)\) & \(1\) & \(\condC{1}\), even  \\
    \(3_{\rm b}\) & \(q\) & \(3\bu_1 \oplus q(-1)\) & \(3\) & -- \\
    \(3_{\rm c}\) & \(q\) & \(2\bu_1 \oplus q(-1)\) & \(3\) & even \\

    \midrule
    
    \(4_{\rm a}\) & \(q\) & \(3\bu_1 \oplus q(-1)\) & \(2\) & \(\condC{2}\), even \\
    \(4_{\rm b}\) & \(q\) & \(2\bu_1 \oplus q(-1)\) & \(4\) & -- \\
    \(4_{\rm c}\) & \(q\) & \(\bu_1 \oplus q(-1)\) & \(4\) & even \\

    \midrule 

    \(5_{\rm a}\) & \(q\) & \(2\bu_1 \oplus q(-1)\) & \(3\) & \(\condC{3}\), even \\
    \(5_{\rm b}\) & \(q\) & \(\bu_1 \oplus q(-1)\) & \(5\) & --\\
    \(5_{\rm c}\) & \(q\) & \(q(-1)\) & \(5\) & even \\

    \midrule
    
    \(6_{\rm a}\) & \(q\) & \(\bu_1 \oplus q(-1)\) & \(4\) & \(\condC{4}\), even \\
    \(6_{\rm b}\) & \(q\) & \(q(-1)\) & \(6\) & -- \\
    \(6_{\rm c}\) & \(\bu_1 \oplus q\) & \(q(-1)\) & \(4\) & even \\
    
    \midrule 
    
    \(7_{\rm a}\) & \(q\) & \(q(-1)\) & \(5\) & \(\condB{5,0}\), \(\condC{5}\), even \\
    \(7_{\rm b}\) & \(\bu_1 \oplus q\) & \(q(-1)\) & \(5\) & \(\condB{5,0}\) \\
    \(7_{\rm c}\) & \(2\bu_1 \oplus q\) & \(q(-1)\) & \(3\) & \(\condB{5,0}\), even \\
    
    \midrule 
    
    \(8_{\rm a}\) & \(\bu_1 \oplus q\) & \(q(-1)\) & \(4\) & \(\condB{4,0}\), \(\condC{4}\), even \\
    \(8_{\rm b}\) & \(2\bu_1 \oplus q\) & \(q(-1)\) & \(4\) & \(\condB{4,0}\) \\
    \(8_{\rm c}\) & \(3\bu_1 \oplus q\) & \(q(-1)\) & \(2\) & \(\condB{4,0}\), even \\
    
    \midrule

    \(9_{\rm a}\) & \(2\bu_1 \oplus q\) & \(q(-1)\) & \(3\) & \(\condB{3,0}\), \(\condC{3}\), even \\
    \(9_{\rm b}\) & \(3\bu_1 \oplus q\) & \(q(-1)\) & \(3\) & \(\condB{3,0}\) \\
    \(9_{\rm c}\) & \(4\bu_1 \oplus q\) & \(q(-1)\) & \(1\) & \(\condB{3,0}\), even \\
    
    \midrule
    
    \(10_{\rm a}\) & \(3\bu_1 \oplus q\) & \(q(-1)\) & \(2\) & \(\condB{2,0}\), \(\condC{2}\), even \\
    \(10_{\rm b}\) & \(4\bu_1 \oplus q\) & \(q(-1)\) & \(2\) & \(\condB{2,0}\) \\
    \(10_{\rm c}\) & \(5\bu_1 \oplus q\) & \(q(-1)\) & \(0\) & \(\condB{2,0}\) \\ 
    
    \midrule

    \(11_{\rm a}\) & \(4\bu_1 \oplus q\) & \(q(-1)\) & \(1\) & \(\condB{1,0}\), \(\condC{1}\), even \\
    \(11_{\rm b}\) & \(5\bu_1 \oplus q\) & \(q(-1)\) & \(1\) & \(\condB{1,0}\) \\  
    
    \midrule 
    
    \(12\) & \(5\bu_1\) & -- & -- & -- \\

    \bottomrule
    \end{tabular}
\end{table}
\end{center}

\subsection{Proof of \texorpdfstring{\Cref{thm:transcendental}}{}} \label{subsec:proof-transc-thm}

The following lemma provides a justification for the terms ``odd'' and ``even'' as applied to quadratic forms.

\begin{lemma}[cf. {\cite[end of p. 130]{Nikulin:int.sym.bilinear.forms}}] \label{lem:l2=rank}
A lattice \(L\) satisfies \(\ell_2(q(L)) = \rank(L)\) if and only if there exists a lattice \(L'\) such that \(L = L'(2)\). Furthermore, \(L'\) is even if and only if \(q(L)\) is even.
\end{lemma}

\begin{proof}
The existence of a lattice \(L'\) with \(L = L'(2)\) is equivalent to the condition \(b(v, w) \in 2\IZ\) for all \(v, w \in L\). This is further equivalent to requiring that \(\frac{1}{2}v \in L^\vee\) for every \(v \in L\), which holds if and only if \(\ell_2(q(L)) = \rank(L)\).

Assume this condition is satisfied. If \(q(L)\) is odd, then there exists an element \(\alpha \in L^\vee/L\) of order \(2\) such that \(q(\alpha) = \frac{1}{2}\varepsilon\), where \(\varepsilon \in \{1, 3\}\). Writing \(\alpha = \frac{1}{2}v\) for some \(v \in L\), we find that \(b(v, v) \equiv 4q(\alpha) = 2\varepsilon \pmod{8}\), implying that \(L'\) is odd.

Conversely, if \(L'\) is odd, then there exists \(v \in L\) such that \(b(v, v) \not\in 4\IZ\). In this case, \(\frac{1}{2}v \in L^\vee/L\) generates a copy of \(\bw_{2,1}^\varepsilon\) in \(q(L)\), confirming that \(q(L)\) is odd.
\end{proof}

Let \(X\) be a K3 surface with transcendental lattice \(T\) of rank~\(\lambda\).
By Keum's criterion (\Cref{thm:keum}), if \(T\) is not a co-idoneal lattice, then \(X\) covers an Enriques surface if and only if there exists a primitive embedding \(T \hookrightarrow \Lminus\).
\Cref{thm:transcendental} follows from \Cref{thm:co-idoneal} once we prove that the conditions given by \Cref{prop:T-embedding-in-Lminus} are equivalent to conditions \ref{cond:2<=lambda<=6}--\ref{cond:lambda=12}.
Therefore, we need to analyze all cases of \Cref{tab:T-embedding-in-Lminus}.

Let us first consider the case \(2 \leq \lambda \leq 6\). We want to prove that \ref{cond:2<=lambda<=6} holds if and only if one of the following holds
\begin{enumerate}[(a)]
 \item \label{cond:(a)} \(\ell_2(q(T)) = \lambda - 2\), \(q(T)\) is even and satisfies condition \(\condC{\lambda - 2}\) (case~\(\lambda_{\rm a})\);
 \item \label{cond:(b)} \(\ell_2(q(T)) = \lambda\) (case \(\lambda_{\rm b}\) or \(\lambda_{\rm c}\)).
\end{enumerate}

Let \(e_1,\ldots,e_ \lambda\) be a system of generators of \(T\).
Suppose first that the corresponding Gram matrix satisfies~\ref{cond:2<=lambda<=6}.

If \(a_{1j}\) is even for \(2 \leq j \leq \lambda\), then \(\ell_2(q(T)) = \rank T = \lambda\) by \Cref{lem:l2=rank}, hence \ref{cond:(b)} holds.

If this is not true, we can suppose \(a_{12}\) to be odd and \(a_{1j}\) to be even for \(3 \leq j \leq \lambda\), up to relabelling and substituting \(e_j\) with \(e_j + e_2\).
Let \(T'\) be the sublattice generated by \(e_1' = 2e_1,e_2,\ldots,e_ \lambda\).
Then \(q(T') \cong \bu_1 \oplus q(T)\), where the copy of \(\bu_1\) is generated by \(e_1'/2\) and \(e_2/2\).
Since \(T' \cong T''(2)\) with \(T''\) even, \(q(T')\) is even and \(\ell_2(q(T')) = \lambda\), by \Cref{lem:l2=rank}. 
Moreover, \(q(T')\) satisfies condition \(\condC{\lambda}\) by \Cref{thm:Nikulin-existence}. 
This implies \ref{cond:(a)}.

Conversely, if \ref{cond:(b)} holds, then \Cref{lem:l2=rank} implies that \(T \cong T'(2)\) for some lattice~\(T'\). 
If \(T'\) is even, then \ref{cond:2<=lambda<=6} holds.
If \(T'\) is odd, then up to relabelling we can suppose that \(e_1^2 \equiv 2 \mod 4\). Then, up to substituting \(e_j\) with \(e_j + e_1\), we can suppose that \(e_j^2 \equiv 0 \mod 4\). Hence, \ref{cond:2<=lambda<=6} holds.

Finally, suppose that \ref{cond:(a)} holds. 
Since \(T\) exists, \(q\) satisfies also conditions \(\condA{4 - \lambda}\) and \(\condB{2, \lambda - 2}\). Therefore, \(\bu_1 \oplus q\) satisfies conditions \(\condA{4 - \lambda}\), \(\condB{2, \lambda - 2}\) and \(\condC{\lambda}\), so the genus~\(g\) of even lattices of signature \((2, \lambda - 2)\) and discriminant form \(\bu_1\oplus q\) is nonempty.
All lattices \(T'\) in \(g\) are of the form \(T' \cong T''(2)\) for some even lattice~\(T''\).
By Proposition~1.4.1 in \cite{Nikulin:int.sym.bilinear.forms}, \(T\) must be an overlattice of such a lattice \(T'\). 
The fact that \(\det T' = 4 \det T\) implies that \(T'\) has index~\(2\) in \(T\). 
Therefore we can find a basis \(e_1, \ldots e_ \lambda\) with \(e_1 \notin T'\) and \(e_j \in T'\) for \(j = 2,\ldots, \lambda\), whose corresponding Gram matrix satisfies~\ref{cond:2<=lambda<=6}.

We now turn to \(\lambda \geq 7\).
The case \(\lambda = 12\) follows immediately from Keum's criterion.
The arguments for \(\lambda \in \{7,\ldots,11\}\) are very similar, so we illustrate here only the case \(\lambda = 10\).

Suppose case \(\lambda_{\rm a}\) holds, i.e., \(\lambda = 10\) and \(q(T) \cong 3\bu_1 \oplus q\), with \(\ell_2(q) = 2\), \(q\) even and satisfying conditions \(\condB{2, 0}\), \(\condC{2}\). 
Since \(T\) exists, \(q\) satisfies also \(\condA{2}\), by \Cref{thm:Nikulin-existence}. 
Hence, using \Cref{thm:Nikulin-existence} again and \Cref{lem:l2=rank}, we infer that there exists an even lattice \(T'\) of signature \((2,0)\) such that \(T'(2)\) has discriminant form~\(q\). 
Since \(T\) is unique in its genus (\Cref{thm:Nikulin-uniqueness}), \(T \cong \Etilde(-1) \oplus T'(2)\), so~\ref{cond:lambda=10.1} holds.

Suppose case \(\lambda_{\rm b}\) or \(\lambda_{\rm c}\) holds, i.e., \(\lambda = 10\) and \(q(T) \cong 4\bu_1 \oplus q\), with \(\ell_2(q) = 2\), and \(q\) satisfying condition \(\condB{2, 0}\).
Since \(T\) exists, \(q\) satisfies also \(\condA{2}\) and \(\condC{2}\). 
Hence, there exists a lattice \(T'\) of signature \((2,0)\) such that \(T'(2)\) has discriminant form \(q\). 
Again by uniqueness, \(T \cong \bE_8(-2) \oplus T'(2)\), so~\ref{cond:lambda=10.2} holds.

Conversely, if \ref{cond:lambda=10.1} holds, then \(q(T) \cong 3\bu_1 \oplus q\), with \(q = q(T')\) being an even finite quadratic form satisfying conditions \(\condA{2}\), \(\condB{2,0}\) and \(\condC{2}\) by \Cref{thm:Nikulin-existence}. 
Hence, by \Cref{thm:Nikulin-existence} and \Cref{prop:Nikulin1.15.1-discriminant-form-orthogonal}, there exists a primitive embedding \(T \hookrightarrow \Lminus\). 
An analogous argument works for \ref{cond:lambda=10.2}. \qed

\begin{remark}
For the equivalence between \ref{cond:lambda=9a} and case \(9_{\rm a}\), one uses the following fact, which is a consequence of \cite[Prop. 1.8.2]{Nikulin:int.sym.bilinear.forms}: if \(q,q'\) are two torsion quadratic forms, then \(\bu_1 \oplus q \cong \bu_1 \oplus q'\) if and only if \(q \cong q'\). 
\end{remark}

\subsection{Co-idoneal lattices} \label{subsec:co-idoneal}
Here, we prove \Cref{thm:co-idoneal} and provide an algorithm to enumerate all co-idoneal lattices (\Cref{def:co-idoneal}).

\begin{lemma} \label{lem:case_lambda_b} 
If \(T\) is a co-idoneal lattice, of rank \(\lambda\), then case \(\lambda_{\rm b}\) of \Cref{tab:cases} holds and \(q(T)\) is odd.
\end{lemma} 
\proof
Consider a primitive embedding \(T \hookrightarrow \Lminus\).
By inspection of \Cref{tab:T-embedding-in-Lminus}, one sees that \(q(T)\) is even if and only if \(q(T^\perp)\) is even. 
In case \(\lambda_{\rm a}\), \(\ell_2(q(T^\perp)) = \rank T^\perp\) and \(q(T)\) is even, so \(T^\perp \cong T'(2)\) for some even lattice \(T'\), by \Cref{lem:l2=rank}. 
Thus, \(T^\perp\) does not contain a vector of square \(-2\) and \(T\) cannot be co-idoneal.
If an embedding as in case~\(\lambda_{\rm c}\) exists, then \(q(T)\) is even and an embedding of \(T\) as in case \(\lambda_{\rm b}\) exists: it suffices to choose a smaller subgroup~\(H\) in \Cref{prop:Nikulin1.15.1-discriminant-form-orthogonal}. Indeed, this changes the type~\((\lambda_{\rm c}, q)\) in \Cref{tab:T-embedding-in-Lminus} to \((\lambda_{\rm b}, \bu_1 \oplus q)\) which does not affect condition \(\condB{12-\lambda_b,0}\). 
Moreover, if an embedding as in case \(\lambda_{\rm b}\) exists and \(q(T)\) is even, then we can argue as before and \(T\) cannot be co-idoneal.
\endproof

Since in case \(\lambda_{\rm b}\) it holds that \(\ell_2(q(T)) = \lambda = \rank T\), \Cref{lem:l2=rank} implies the following corollary.

\begin{corollary} \label{cor:half.exc} 
If \(T\) is a co-idoneal lattice, then \(T \cong T'(2)\) for some odd lattice~\(T'\) (which we call \emph{co-idoneal lattice half}). \qed
\end{corollary} 

The following proposition explains the connection between co-idoneal lattices and idoneal genera.

\begin{proposition} \label{prop:co-idoneal.lattices/idoneal.genera}
For each co-idoneal lattice \(T\) there exists a unique idoneal genus~\(g\) with the following property: for each primitive embedding \(T \hookrightarrow \Lminus\), there exists a lattice \(L\) in \(g\) with \(T^\perp \cong L(-2)\).
\end{proposition}
\proof
Let \(\lambda = \rank T\) and consider a primitive embedding \(T \hookrightarrow \Lminus\). 
It follows from \Cref{lem:case_lambda_b} that \(\ell_2(q(T^\perp)) = 12 - \lambda = \rank T^\perp\) and \(q(T)\) is odd. 
By \Cref{lem:l2=rank}, there exists an odd lattice \(L\) with \(T^\perp = L(-2)\). Let \(g\) be the genus of \(L\).

The discriminant form of \(T^\perp\) is determined by the discriminant form of~\(T\) according to \Cref{tab:T-embedding-in-Lminus}. 
According to \cite[Corollary~1.16.3]{Nikulin:int.sym.bilinear.forms}, the genus of a lattice is determined by its signature, parity and discriminant bilinear form. Hence, each lattice \(L\) such that \(T^\perp \cong L(-2)\) belongs to the same genus~\(g\). 

Conversely, as \(\Lminus\) is unique in its genus, each lattice \(L\) in \(g\) satisfies \(L(-2) \cong T^\perp\) for some embedding \(T \hookrightarrow \Lminus\), by \Cref{prop:Nikulin1.15.1-discriminant-form-orthogonal}.
Since \(T^\perp\) always contains a vector of square \(-2\), each lattice \(L\) in \(g\) contains a vector of square \(1\), i.e., \(g\) is idoneal.
\endproof

\proof[Proof of \Cref{thm:co-idoneal}]
By \Cref{thm:idoneal}, there exists only finitely many idoneal genera. Each idoneal genus \(g\) consists of finitely many isomorphim classes of lattices~\(L\). For each such lattice~\(L\), there exists finitely isomorphism classes of lattices~\(T\) such that \(L(-2)^\perp \cong T\) for some primititive embedding \(L(-2) \hookrightarrow \Lminus\). By \Cref{prop:co-idoneal.lattices/idoneal.genera}, we conclude that there exist finitely many co-idoneal lattices.
\endproof 

\begin{theorem} \label{thm:algorithm.co-idoneal}
    For all \(2 \leq \lambda \leq 11\), there exists an algorithm to compute the list of all co-idoneal lattices \(T\) of rank \(\lambda\).
\end{theorem}

\begin{proof}
Consider the following algorithm:
    \begin{myalgorithm}\caption{Co-idoneal lattices}\label{alg:co-idoneal}
    \begin{algorithmic}[1]
         \REQUIRE an integer \(\lambda\) with \(2 \leq \lambda \leq 11\)
         \ENSURE the list of all co-idoneal lattice halves of rank \(\lambda\)

         \STATE \label{line:idoneal} \(I \longleftarrow\) list of idoneal genera of rank~\(n = 12-\lambda\)
         \STATE \(H \longleftarrow \emptyset\)
         \FORALL{\(g \in I\)}
            \STATE \(L \longleftarrow\) a representative of \(g\)
            \STATE \label{line:x} \(x \longleftarrow q(L(-2))\)
            \IF{\(x\) is one of the suitable forms \(q(T^\perp)\) in \Cref{tab:T-embedding-in-Lminus}}
                \STATE \label{line:y} \(y \longleftarrow\) quadratic form corresponding to \(q(T)\) given by \Cref{tab:T-embedding-in-Lminus}
                \STATE \label{line:C} \(C \longleftarrow\) list of even lattices \(T\) of signature \((2,\lambda -2)\) and discriminant form \(y\)
                \FORALL{\(T \in C\)}
                    \STATE \(T' \longleftarrow T(1/2)\) \label{step:T(1/2)}
                    \STATE \(H \longleftarrow H \cup \{T'\}\)
                \ENDFOR
            \ENDIF
         \ENDFOR 
         \RETURN \(H\)
    \end{algorithmic}
    \end{myalgorithm}
    The list \(I\) in line \ref{line:idoneal} is produced using \Cref{thm:algorithm.idoneal}.
    The algorithm is justified by \Cref{prop:T-embedding-in-Lminus,prop:co-idoneal.lattices/idoneal.genera}.  The fact that in line~\ref{step:T(1/2)} we obtain an integral lattice is justified by \Cref{cor:half.exc}.
\end{proof}

\begin{example}
We illustrate the algorithm from \Cref{thm:algorithm.co-idoneal} for the case \(\lambda = 2\).

In line~\ref{line:idoneal}, the algorithm generates the list of all idoneal genera of rank \(12 - \lambda = 10\). This list, available in the file \verb+idoneal.genera.txt+ \cite{brandhorst_2024_10617125}, comprises exactly six genera, represented by the following lattices:
\begin{align*}
    L_1 &= \bE_8 \oplus 2[1], & 
    L_2 &= \bE_8 \oplus [2] \oplus [1], &
    L_3 &= \bA_2 \oplus 8[1], \\
    L_4 &= \bD_9 \oplus [1], &
    L_5 &= \bA_3 \oplus 7[1], &
    L_6 &= \bD_8 \oplus [2] \oplus [1].
\end{align*}

Next, in line~\ref{line:x}, the algorithm computes the quadratic forms \(x_i = q(L_i(-2))\) for \(i = 1, \ldots, 6\), yielding:
\begin{align*}
    x_1 &= 4\bu_1 \oplus 2\bw_{2,1}^3, & 
    x_2 &= 4\bu_1 \oplus \bw_{2,1}^3 \oplus \bw_{2,2}^7, \\
    x_3 &= 3\bu_1 \oplus \bv_1 \oplus \bw_{2,1}^1 \oplus \bw_{2,1}^3 \oplus \bw_{3,1}^{-1}, &
    x_4 &= 4\bu_1 \oplus \bw_{2,1}^3 \oplus \bw_{2,3}^7, \\
    x_5 &= 3\bu_1 \oplus \bv_1 \oplus \bw_{2,1}^1 \oplus \bw_{2,3}^1, &
    x_6 &= 2\bu_1 \oplus \bv_1 \oplus \bw_{2,1}^1 \oplus \bu_2 \oplus \bw_{2,2}^1.
\end{align*}

Among these, only \(x_1\), \(x_2\), and \(x_4\) are of the form \(q(T^\perp)\) described in \Cref{tab:T-embedding-in-Lminus}. Specifically, they take the form \(4\bu_1 \oplus q_i(-1)\) of case~\(2_\mathrm{b}\). In line~\ref{line:y} of \Cref{alg:co-idoneal}, the corresponding forms \(y_i\) are determined from the column \(q(T)\) in \Cref{tab:T-embedding-in-Lminus}. Thus, \(y_i = q_i\) and we have
\[
    y_1 = 2\bw_{2,1}^1, \qquad 
    y_2 = \bw_{2,1}^1 \oplus \bw_{2,2}^1, \qquad  
    y_4 = \bw_{2,1}^1 \oplus \bw_{2,3}^1.
\]

In line~\ref{line:C}, the algorithm computes all even lattices of signature \((2, 0)\) and discriminant forms \(y_1, y_2, y_4\). Each lattice is unique in its genus, given by:
\[
    T_1 = 2[2], \qquad T_2 = [2] \oplus [4], \qquad T_4 = [2] \oplus [8].
\]

Thus, these are all co-idoneal lattices of rank \(\lambda = 2\). Their halves are documented in the file \verb+half.co-idoneal.lat.txt+ \cite{brandhorst_2024_10617125}. In this way, we retrieve the original result of Sertöz~\cite{Sertoez05}.
\end{example}

At the end of our case-by-case analysis, the following fact turns out to be true.

\begin{addendum}
All co-idoneal lattices are unique in their genus.
\qed
\end{addendum}

\printbibliography

@article {Allcock.Gal.Mark,
    AUTHOR = {Allcock, Daniel and Gal, Itamar and Mark, Alice},
     TITLE = {The {C}onway-{S}loane calculus for 2-adic lattices},
   JOURNAL = {Enseign. Math.},
  FJOURNAL = {L'Enseignement Math\'{e}matique},
    VOLUME = {66},
      YEAR = {2020},
    NUMBER = {1-2},
     PAGES = {5--31},
      %ISSN = {0013-8584,2309-4672},
   MRCLASS = {11E08},
  MRNUMBER = {4162280},
MRREVIEWER = {Jingbo\ Liu},
       DOI = {10.4171/LEM/66-1/2-2},
%       URL = {https://doi.org/10.4171/LEM/66-1/2-2},
}

@phdthesis{kirschmer,
  address = {Aachen, Germany},
  author = {Kirschmer, Markus},
  school = {RWTH Aachen University},
  title = {Definite quadratic and hermitian forms with small class number},
  type = {Habilitationsschrift},
  year = {2016},
  URL = {http://www.math.rwth-aachen.de/~Markus.Kirschmer/papers/herm.pdf},
}

@book {landau,
    AUTHOR = {Landau, Edmund},
     TITLE = {Handbuch der {L}ehre von der {V}erteilung der {P}rimzahlen. 2
              {B}\"{a}nde},
      NOTE = {2d ed,
              with an appendix by Paul T. Bateman},
 PUBLISHER = {Chelsea Publishing Co., New York},
      YEAR = {1953},
     PAGES = {xviii+pp. 1--564; ix+pp. 565--1001},
   MRCLASS = {10.0X},
  MRNUMBER = {68565},
MRREVIEWER = {L.\ Schoenfeld},
}

@article {ConwaySloane88,
    AUTHOR = {Conway, John H. and Sloane, Neil J. A.},
     TITLE = {Low-dimensional lattices. {IV}. {T}he mass formula},
   JOURNAL = {Proc. Roy. Soc. London Ser. A},
  FJOURNAL = {Proceedings of the Royal Society. London. Series A.
              Mathematical, Physical and Engineering Sciences},
    VOLUME = {419},
      YEAR = {1988},
    NUMBER = {1857},
     PAGES = {259--286},
%      ISSN = {0962-8444},
   MRCLASS = {11H55 (11E32 11H06)},
  MRNUMBER = {965484},
MRREVIEWER = {O. H. K\"{o}rner},
    url = {https://www.jstor.org/stable/2398465},
}

@article {cho2015,
    AUTHOR = {Cho, Sungmun},
     TITLE = {Group schemes and local densities of quadratic lattices in
              residue characteristic 2},
   JOURNAL = {Compos. Math.},
  FJOURNAL = {Compositio Mathematica},
    VOLUME = {151},
      YEAR = {2015},
    NUMBER = {5},
     PAGES = {793--827},
      %ISSN = {0010-437X,1570-5846},
   MRCLASS = {11E41 (11E12 11E57 11E95 14L15 20G25)},
  MRNUMBER = {3347991},
MRREVIEWER = {Hidenori\ Katsurada},
       DOI = {10.1112/S0010437X14007829},
       %URL = {https://doi.org/10.1112/S0010437X14007829},
}

@article{Johansson2017arb,
  author = {Johansson, Fredrik},
  title = {Arb: efficient arbitrary-precision midpoint-radius interval arithmetic},
  journal = {IEEE Transactions on Computers},
  year = {2017},
  volume = {66},
  issue = {8},
  pages = {1281--1292},
  doi = {10.1109/TC.2017.2690633},
}

@Book{splag,
 Author = {J. H. {Conway} and N. J. A. {Sloane}},
 Title = {{Sphere packings, lattices and groups}},
 FJournal = {{Grundlehren der Mathematischen Wissenschaften}},
 Journal = {{Grundlehren Math. Wiss.}},
 ISSN = {0072-7830},
 Volume = {290},
 ISBN = {0-387-98585-9},
 Pages = {lxxiv, 703},
 Year = {1999},
 Publisher = {New York, NY, Springer},
 MSC2010 = {52C17 52C07 11H06 11-02 05-02 05B40 52-02 03B30 05B05 20D08 11F27 94B99},
 Zbl = {0915.52003}
}

@article {euler:numeri.idonei,
    author =        "Euler, Leonhard",
    title =         "Illustratio paradoxi circa progressionem numerorum idoneorum sive congruorum",
    journal =       "Nova Acta Academiae Scientarum Imperialis Petropolitinae",
    volume =        "15",
    year =          "1806", 
    pages =         "29--32",
    language =      "Latin",
    note =          "Delivered to the St. Petersburg Academy April 20, 1778. English translation by Jordan Bell: \href{https://arxiv.org/abs/math/0507352}{\texttt{arXiv:math/0507352}} (2005)",
}

@article {Festi.Veniani,
    AUTHOR = {Festi, Dino and Veniani, Davide Cesare},
     TITLE = {Enriques involutions on pencils of {K}3 surfaces},
   JOURNAL = {Math. Nachr.},
  FJOURNAL = {Mathematische Nachrichten},
    VOLUME = {295},
      YEAR = {2022},
    NUMBER = {7},
     PAGES = {1312--1326},
%      ISSN = {0025-584X,1522-2616},
   MRCLASS = {14J28 (14J27 14J50)},
  MRNUMBER = {4468868},
MRREVIEWER = {Remke\ Kloosterman},
DOI = {doi.org/10.1002/mana.202100140},
%URL = {https://doi.org/10.1002/mana.202100140},
}

@book {hua,
    AUTHOR = {Hua, Loo Keng},
     TITLE = {Introduction to number theory},
      NOTE = {Translated from the Chinese by Peter Shiu},
 PUBLISHER = {Springer-Verlag, Berlin-New York},
      YEAR = {1982},
     PAGES = {xviii+572},
      %ISBN = {3-540-10818-1},
   MRCLASS = {10-01},
  MRNUMBER = {665428},
}

@article {Kani:idoneal.numbers,
    AUTHOR = {Kani, Ernst},
     TITLE = {Idoneal numbers and some generalizations},
   JOURNAL = {Ann. Sci. Math. Qu\'{e}bec},
  FJOURNAL = {Annales des Sciences Math\'{e}matiques du Qu\'{e}bec},
    VOLUME = {35},
      YEAR = {2011},
    NUMBER = {2},
     PAGES = {197--227},
      %ISSN = {0707-9109},
   MRCLASS = {11R11 (01A50 11-03 11E41 11R29)},
  MRNUMBER = {2917832},
MRREVIEWER = {St\'{e}phane R. Louboutin},
       URL = {https://www.labmath.uqam.ca/~annales/volumes/35-2/PDF/197-227.pdf},
}

@article {Keum90, 
    AUTHOR = {Keum, JongHae},
     TITLE = {Every algebraic {K}ummer surface is the {$K3$}-cover of an
              {E}nriques surface},
   JOURNAL = {Nagoya Math. J.},
  FJOURNAL = {Nagoya Mathematical Journal},
    VOLUME = {118},
      YEAR = {1990},
     PAGES = {99--110},
      %ISSN = {0027-7630},
   MRCLASS = {14J28 (14J50)},
  MRNUMBER = {1060704},
MRREVIEWER = {Shigeyuki Kondo},
       DOI = {10.1017/S0027763000003019},
%       URL = {https://doi.org/10.1017/S0027763000003019},
}

@book {kneser,
    AUTHOR = {Kneser, Martin},
     TITLE = {Quadratische {F}ormen},
      NOTE = {Revised and edited in collaboration with Rudolf Scharlau},
 PUBLISHER = {Springer-Verlag, Berlin},
      YEAR = {2002},
     PAGES = {viii+164},
      %ISBN = {3-540-64650-7},
   MRCLASS = {11-02 (11Exx 11H55)},
  MRNUMBER = {2788987},
       DOI = {10.1007/978-3-642-56380-5},
%       URL = {https://doi.org/10.1007/978-3-642-56380-5},
    language =  "German",       
}

@article {Lee12,
    AUTHOR = {Lee, Kwangwoo},
     TITLE = {Which {K}3 surfaces with {P}icard number 19 cover an
              {E}nriques surface},
   JOURNAL = {Bull. Korean Math. Soc.},
  FJOURNAL = {Bulletin of the Korean Mathematical Society},
    VOLUME = {49},
      YEAR = {2012},
    NUMBER = {1},
     PAGES = {213--222},
      %ISSN = {1015-8634},
   MRCLASS = {14J28},
  MRNUMBER = {2931561},
MRREVIEWER = {Hisanori Ohashi},
       DOI = {10.4134/BKMS.2012.49.1.213},
%       URL = {https://doi.org/10.4134/BKMS.2012.49.1.213},
}

@article {morrison:K3.large.picard,
    AUTHOR = {Morrison, David R.},
     TITLE = {On {$K3$} surfaces with large {P}icard number},
   JOURNAL = {Invent. Math.},
  FJOURNAL = {Inventiones Mathematicae},
    VOLUME = {75},
      YEAR = {1984},
    NUMBER = {1},
     PAGES = {105--121},
      %ISSN = {0020-9910},
   MRCLASS = {14J28 (14C30 14J05 14K10)},
  MRNUMBER = {728142},
MRREVIEWER = {I. Dolgachev},
       DOI = {10.1007/BF01403093},
%       URL = {https://doi.org/10.1007/BF01403093},
}

@article{nikulin:kummer,
    author = {Nikulin, Viacheslav V.},
    title = {On Kummer surfaces},
    journal = {Math. USSR-Izv.},
    volume = {9},
    number = {2},
    year = {1975},
    pages = {261--275},
    DOI = {10.1070/IM1975v009n02ABEH001477},
}

@article{Nikulin:int.sym.bilinear.forms,
    author = {Nikulin, Viacheslav V.},
    title = {Integral symmetric bilinear forms and some of their applications},
    journal = {Math. USSR-Izv.},
    year = {1980},
    volume = {14},
    number = {1},
    pages = {103--167},
    doi = {10.1070/IM1980v014n01ABEH001060},
}

@article {Nikulin:int.sym.bilinear.forms_Russian,
    AUTHOR = {Nikulin, Viacheslav V.},
     TITLE = {Integer symmetric bilinear forms and some of their geometric
              applications},
   JOURNAL = {Izv. Akad. Nauk SSSR Ser. Mat.},
  FJOURNAL = {Izvestiya Akademii Nauk SSSR. Seriya Matematicheskaya},
    VOLUME = {43},
      YEAR = {1979},
    NUMBER = {1},
     PAGES = {111--177},
     % ISSN = {0373-2436},
   MRCLASS = {10C05 (14G30 14J17 14J25 57M99 57R45 58C27)},
  MRNUMBER = {525944},
MRREVIEWER = {I. Dolgachev},
  language = "Russian",
    doi = {10.1070/IM1980v014n01ABEH001060},
}

@article {neighbor,
    AUTHOR = {Scharlau, Rudolf and Hemkemeier, Boris},
     TITLE = {Classification of integral lattices with large class number},
   JOURNAL = {Math. Comp.},
  FJOURNAL = {Mathematics of Computation},
    VOLUME = {67},
      YEAR = {1998},
    NUMBER = {222},
     PAGES = {737--749},
      %ISSN = {0025-5718},
   MRCLASS = {11E41 (11H55)},
  MRNUMBER = {1458224},
MRREVIEWER = {Renaud Coulangeon},
       DOI = {10.1090/S0025-5718-98-00938-7},
       %URL = {https://doi.org/10.1090/S0025-5718-98-00938-7},
}

@misc{hanke,
AUTHOR = {Hanke, Jonathan and {The Sage Development Team}},
TITLE = "Quadratic Forms in {SageMath}",
note = "\url{https://doc.sagemath.org/html/en/reference/quadratic_forms/sage/quadratic_forms/quadratic_form.html}",
YEAR = {2007},
}

@misc{brandhorst_discriminant_forms,
AUTHOR = {Brandhorst, Simon and {The Sage Development Team}},
TITLE = "Quadratic lattices and discriminant forms in {SageMath}",
note = "\url{https://doc.sagemath.org/html/en/reference/modules/sage/modules/torsion_quadratic_module.html}, \newline \url{https://doc.sagemath.org/html/en/reference/modules/sage/modules/free_quadratic_module_integer_symmetric.html}" ,
YEAR = {2017},
}

@article{fincke.pohst,
    AUTHOR = {Fincke, U. and Pohst, M.},
     TITLE = {Improved methods for calculating vectors of short length in a
              lattice, including a complexity analysis},
   JOURNAL = {Math. Comp.},
  FJOURNAL = {Mathematics of Computation},
    VOLUME = {44},
      YEAR = {1985},
    NUMBER = {170},
     PAGES = {463--471},
      %ISSN = {0025-5718},
   MRCLASS = {11H50 (11H55)},
  MRNUMBER = {777278},
MRREVIEWER = {A. K. Lenstra},
       DOI = {10.2307/2007966},
       %URL = {https://doi.org/10.2307/2007966},
}

@misc{OEIS:A000926,
    author =    "Sloane, Neil J. A.",
    title =     "Euler's idoneal numbers, {Entry} \href{http://oeis.org/A000926}{\texttt{A000926}} in {The On-Line Encyclopedia of Integer Sequences}",
    note =  "\href{http://oeis.org/A000926}{\texttt{oeis.org/A000926}}",
    year =      "2020",
}

@article {isometries,
    AUTHOR = {Plesken, Wilhelm and Souvignier, Bernd},
     TITLE = {Computing isometries of lattices},
      NOTE = {Computational algebra and number theory (London, 1993)},
   JOURNAL = {J. Symbolic Comput.},
  FJOURNAL = {Journal of Symbolic Computation},
    VOLUME = {24},
      YEAR = {1997},
    NUMBER = {3-4},
     PAGES = {327--334},
      %ISSN = {0747-7171},
   MRCLASS = {11H56 (11Y16)},
  MRNUMBER = {1484483},
MRREVIEWER = {Christine Bachoc},
       DOI = {10.1006/jsco.1996.0130},
%       URL = {https://doi.org/10.1006/jsco.1996.0130},
}

@article {magma,
    AUTHOR = {Bosma, Wieb and Cannon, John and Playoust, Catherine},
     TITLE = {The {M}agma algebra system. {I}. {T}he user language},
      NOTE = {Computational algebra and number theory (London, 1993)},
   JOURNAL = {J. Symbolic Comput.},
  FJOURNAL = {Journal of Symbolic Computation},
    VOLUME = {24},
      YEAR = {1997},
    NUMBER = {3-4},
     PAGES = {235--265},
      %ISSN = {0747-7171},
   MRCLASS = {68Q40},
       DOI = {10.1006/jsco.1996.0125},
%       URL = {http://dx.doi.org/10.1006/jsco.1996.0125},
}

@manual{sage,
  Key          = {Sage},
  Author       = {Stein, William A. and others},
  Organization = {The Sage Development Team},
  Title        = {{S}age {M}athematics {S}oftware ({V}ersion 9.0)},
  note         = {\url{www.sagemath.org}},
  Year         = {2020},
}

@manual{pari,
    organization = "{The PARI~Group}",
    title        = "{PARI/GP version {\tt 2.11.1}}",
    year         = "2018",
    address      = "Univ. Bordeaux",
    note         = "\url{pari.math.u-bordeaux.fr/}"
    }

@article {Ohashi07,
    AUTHOR = {Ohashi, Hisanori},
     TITLE = {On the number of {E}nriques quotients of a {$K3$} surface},
   JOURNAL = {Publ. Res. Inst. Math. Sci.},
  FJOURNAL = {Kyoto University. Research Institute for Mathematical
              Sciences. Publications},
    VOLUME = {43},
      YEAR = {2007},
    NUMBER = {1},
     PAGES = {181--200},
%      ISSN = {0034-5318},
   MRCLASS = {14J28},
  MRNUMBER = {2319542},
MRREVIEWER = {Paolo Stellari},
       URL = {https://www.kurims.kyoto-u.ac.jp/~prims/pdf/43-1/43-1-10.pdf},
}

@article {ShimadaVeniani,
    AUTHOR = {Shimada, Ichiro and Veniani, Davide Cesare},
     TITLE = {Enriques involutions on singular {K}3 surfaces of small
              discriminants},
   JOURNAL = {Ann. Sc. Norm. Super. Pisa Cl. Sci. (5)},
  FJOURNAL = {Annali della Scuola Normale Superiore di Pisa. Classe di
              Scienze. Serie V},
    VOLUME = {21},
      YEAR = {2020},
     PAGES = {1667--1701},
%      ISSN = {0391-173X,2036-2145},
   MRCLASS = {14J28 (14J50)},
  MRNUMBER = {4288644},
MRREVIEWER = {Shigeyuki\ Kondo},
       DOI = {10.2422/2036-2145.201902\_004},
%       URL = {https://doi.org/10.2422/2036-2145.201902_004},
}

@article {Ohashi,
    author =        "Ohashi, Hisanori",
    title =         "Counting {E}nriques quotients of a {K3} surface",
    journal =       "Res. Inst. Math. Sci. preprint",
    number =        "1609", 
    note =          "\url{www.kurims.kyoto-u.ac.jp/preprint/file/RIMS1609.pdf}",
    year =          "2007",
    pages =         "1--6",
}

@article {Sertoez05,
    AUTHOR = {Sert\"{o}z, Ali Sinan},
     TITLE = {Which singular {$K3$} surfaces cover an {E}nriques surface},
   JOURNAL = {Proc. Amer. Math. Soc.},
  FJOURNAL = {Proceedings of the American Mathematical Society},
    VOLUME = {133},
      YEAR = {2005},
    NUMBER = {1},
     PAGES = {43--50},
%      ISSN = {0002-9939},
   MRCLASS = {14J28 (11E39)},
  MRNUMBER = {2085151},
MRREVIEWER = {Jong Hae Keum},
       DOI = {10.1090/S0002-9939-04-07666-X},
%       URL = {https://doi.org/10.1090/S0002-9939-04-07666-X},
}

@article {weinberger,
    AUTHOR = {Weinberger, Peter J.},
     TITLE = {Exponents of the class groups of complex quadratic fields},
   JOURNAL = {Acta Arith.},
  FJOURNAL = {Polska Akademia Nauk. Instytut Matematyczny. Acta Arithmetica},
    VOLUME = {22},
      YEAR = {1973},
     PAGES = {117--124},
%      ISSN = {0065-1036},
   MRCLASS = {12A25},
  MRNUMBER = {313221},
MRREVIEWER = {H. Yokoi},
       DOI = {10.4064/aa-22-2-117-124},
%       URL = {https://doi.org/10.4064/aa-22-2-117-124},
}

@book {Whittaker.Watson:course.modern.analysis,
    AUTHOR = {Whittaker, E. T. and Watson, G. N.},
     TITLE = {A course of modern analysis---an introduction to the general theory of infinite processes and of analytic functions with an account of the principal transcendental functions},
    EDITOR = {Moll, Victor H.},
   EDITION = {5th ed.},
      NOTE = {With a foreword by S. J. Patterson},
 PUBLISHER = {Cambridge University Press, Cambridge},
      YEAR = {2021},
     PAGES = {lii+668},
      %ISBN = {978-1-316-51893-9},
   MRCLASS = {30-01 (01A75 30-03)},
  MRNUMBER = {4286926},
}

@misc{Yoruk,
    author =        "Yörük, Oğuzhan",
    title =         "Which algebraic {K3} surfaces doubly cover an {E}nriques surface: a computational approach",
    year =          "2019",
    note =          "Master thesis, Bilkent University",
    howpublished = "\url{http://repository.bilkent.edu.tr/handle/11693/50635}",
}

@misc{brandhorst_2024_10617125,
  author       = {Brandhorst, Simon and
                  Sonel, Serkan and Veniani, Davide Cesare},
  title        = {{Data for "Idoneal genera and K3 surfaces covering 
                   an Enriques surface"}}, 
month        = mar,
  year         = 2020,
    note = {Zenodo},  
publisher    = {Zenodo},
  doi          = {10.5281/zenodo.10617125},
%  url          = {https://doi.org/10.5281/zenodo.10617125}
}

\end{document}